\magnification=\magstep1
\input amstex
\voffset=-3pc

\documentstyle{amsppt}

\vsize=9.5truein
\hsize=6.5truein

\def\bR{\Bbb R}

\def\cU{\Cal U}

\redefine\sM{\Cal M}

\def\uy{\underline y}

\def\dist{\text{dist}}

\def\val{\text{val}}
\def\sS{\Cal S}
\def\sI{ d}
\def\sG{\Cal G}
\def\wt{\text{wt}}
\def\size{\text{size}}
\def\st{\text{st}}
\def\tes{\text{tes}}
\def\cyc{\text{cyc}}
\def\path{\text{path}}
\def\ev{\text{ev}}
\def\ASC{\text{(ASC)}}
\def\valsm{\text{val}_{sm}}
\def\b0{\bold 0}
\def\ua{\underline a}
\def\ub{\underline b}
\def\ux{\underline x}
\def\sM{\Cal M}
\def\bR{\Bbb R}
\def\bF{\Bbb F}
\def\bn{\bold n}
\def\bg{\bold g}
\def\bj{\bold j}
\def\bk{\bold k}
\def\bh{\bold h}
\def\bi{\bold i}

\def\b1{\bold 1}

\magnification=\magstep1
\parskip=6pt
\NoBlackBoxes
\topmatter
\title Graph Games and the Pizza Problem\endtitle  
\author Daniel E.~Brown and Lawrence G.~Brown
\endauthor
%\endtopmatter

\abstract{We propose a class of two person perfect information games based on weighted graphs.
One of these games can be described in terms of a round pizza which is cut radially into pieces of varying size.
The two players alternately take pieces subject to the following rule:\ \ Once the first piece has been chosen, all subsequent selections must be adjacent to the hole left by the previously taken pieces.
Each player tries to get as much pizza as possible.
The original pizza problem was to settle the conjecture that Player One can always get at least half of the pizza.
The conjecture turned out to be false.
Our main result is a complete analysis of a somewhat simpler class of games, concatenations of stacks and two--ended stacks, and we provide a linear time algorithm for this.
The algorithm and its output can be described without reference to games.
It produces a certain kind of partition of a given finite sequence of real numbers.
The conditions on the partition involve alternating sums of various segments of the given sequence.
We do not know whether these partitions have applications outside of game theory.
The algorithm leads to a quadratic time algorithm which gives the value and an optimal initial move for pizza games.
We also provide some general theory concerning the semigroup of equivalence classes of graph games.}
\endabstract
\endtopmatter


\bigskip\noindent
{\bf \S0. Introduction}.

Graph games are a class of two person, perfect information, zero sum games in which the players alternately remove vertices from a finite weighted graph according to certain selection rules.
The main selection rule is that after a vertex has been taken, its neighbors become eligible for selection on subsequent moves (the full definition will be given in the next section).
The game ends when all vertices have been removed, and both players try to maximize the sum of the weights of the vertices that they take.
The outcome for either player is the sum of the weights of the vertices taken by the player minus the sum of the weights of the vertices taken by the opponent.
This causes the game to be a zero sum game rather than just a constant sum game.

Here are some basic examples which provided our main motivation:

1.\ Multiple stacks:\ \ The graph is a disjoint union of paths, thought of as vertical, and on any turn a player may select the top remaining vertex from any path.

2.\ Two--ended stacks:\ \ The graph is a path, thought of as horizontal, and on any turn a player may select the leftmost or rightmost remaining vertex.
Thus if the path has $n$ vertices, there are $2^{n-1}$ possible plays of the game.

3.\ Cycles:\ \ The graph is a cycle.
For the first move Player One may select any vertex.
Then the game turns into a two--ended stack in accordance with the rules stated above.
If the cycle has $n$ vertices, there are $n2^{n-2}$ possible plays of the game.

The pizza games described in the abstract are cycle games (with positive weights).
One can give food sharing interpretations of stacks and two--ended stacks by thinking of the vertices as cookies of varying size.
The original motivation for the selection rules was that the players want to be outwardly polite while trying to get as much food as possible.

In the case of cycle games it would appear that Player One has a natural advantage because there are $n$ choices for the first move and at most two for any subsequent move.
Also if $n$ is odd and the weights are positive, the fact that Player One takes more than half of the vertices would seem to be an advantage.
This intuition underlay our original conjecture that cycle games with non--negative weights have non--negative value, but the conjecture was false. (By {\it value} we mean the optimal outcome for Player One.  Note that this term is not used the same way as in combinatorial game theory.)
The smallest counterexample has 15 vertices, and one possible weight sequence is $(0,1,0,1,0,0,1,0,2,0,0,2,0,2,0)$ (regarded as cyclic).
The fact that some weights are 0 is unimportant, since the value is a continuous function of the weights.

On the other hand, cycle games with an even number of vertices always do have non--negative value because of the availability of a ``color strategy''.
Assume the vertices are alternately colored red or blue.
Player One can arrange to get either all of the red vertices or all of the blue ones.
Color strategies are also available for two--ended stacks with an even number of vertices, and they are available to Player Two, once the initial move has been made, for odd cycles or odd two--ended stacks.

It is natural to consider greedy strategies.
A greedy strategy is one where a player always chooses a vertex of maximal weight from among those which may legally be chosen.
Greedy strategies are optimal for $n$--cycles with $n\leq 4$ and for two--ended stacks with at most three vertices, but neither color strategies nor greedy strategies are optimal in general.
For two--ended stacks:\ \ If the weight sequence is $(4,3,1,2)$, the greedy strategy is optimal and the color strategy is not; if the weight sequence is $(4,3,6,5)$, the color strategy is optimal and the greedy strategy is not; and if the weight sequence is $(4,3,1,4,7,5)$, neither a greedy nor a color strategy is optimal.
(All three of these games have value 2.)
For a 6--cycle with weight sequence $(0,N,2,1,2,1)$, $N>2$, or $(0,3,3,1,2,1)$ both the greedy initial move and the color strategy are wrong.
(The values of these games are $N$ and 2.)
It is also possible to give examples of two--ended stacks or cycles with five vertices where the greedy strategy is wrong for Player One.
Consider the weight sequences $(2,4,1,0,1)$ and $(2,3,1,2,0)$, respectively, with values 0 and 2.
An interesting example is a $(2k+1)$--cycle with weight sequence $(N,1,0,1,0,\ldots,1,0)$, where $N$ is a large positive integer and $k>N$.
Here the greedy initial move is wrong and the value is 0 or 1, according as $N+k$ is even or odd.

The plan of the paper is as follows.
Section One contains definitions, notations, and other preliminaries.
Section Two contains the solution of concatenations of stacks and two--ended stacks.
Section Three contains some general theory about the semigroup of equivalence classes of graph games.
Although not all classes have inverses, every invertible equivalence class is equal to its own inverse.
Section Four applies the earlier results to cycle games and to a related class called unbroached path games.
And Section Five gives a brief discussion, without proofs, of some variants of graph games.

Our approach to this project was heavily influenced by Andrew Gleason's filmed lecture, [2], on Nim--type games.
Thus [2] defines an equivalence relation for Nim--type games and an addition operation on equivalence classes, and includes the Sprague--Grundy Theorem ([3],[6]).
Our main result in \S2, that every stack or two--ended stack is equivalent to a very simple graph game, is our best attempt at an analogue of the Sprague--Grundy Theorem.

The games proposed in this paper were invented by the first author in 1996 and almost all of our joint research was completed by 1998.  The only important exceptions are 4.7 through 4.10, which were done in 2015.  The other exceptions are Remark 1.2. parts of Example 1.7, Corollary 2.12(ii), Remark 2.13. parts of subsection 3.27, Remark 4.13. and parts of $\S$5.  Some other people were told about these games early on, and some other papers about them have been published in the meantime.  Many of these papers concern either games played on graphs with different rules from our graph games or the question what proportion of a pizza Player One can always be sure to get.  J. Cibulka, J.Kyn\v{c}l, V. M\'{e}sz\'{a}ros, R. Stola\'{r}, and P. Valtr [1] and K. Knauer, P. Micek, and T. Ueckerdt [5] independently answered this question by verifying a conjecture of P. Winkler that the optimal proportion is 4/9.  Although we did not originally work on this proportionality question, we did find in 1998 both the pizza with 15 pieces mentioned above and another one mentioned below in 3.27 which has 21 pieces all of weight 0 or 1.  Both of these happen to be extremal for the proportionality problem, but our motivation was just to find relatively simple examples of pizza games with negative value.  Subsection 4.7 below sketches a new proof of the proportionality result of [1] and [5] and also discusses two alternative pizza-sharing games.  Theorem 4.9 was inspired by a question asked near the end of [5].

\bigskip\noindent
{\bf \S1.\ Definitions and basic properties}.

A {\it graph game board} is a pair $g=(G,A)$, where $G$ is a finite vertex--weighted graph and $A$ is a set of vertices of $G$.
Either $A$ or $G$ may be empty and the weight of a vertex $v$, denoted by $\wt(v)$, may be any real number.
The vertices in $A$ are called {\it available}.
For a set $S$ of vertices of $G$, $g\backslash S=(\tilde G,\tilde A)$, where $\tilde G$ is obtained from $G$ be removing the vertices in $S$ and all edges incident to them and 
$\tilde A=(A\cup \{v\colon v$ is a vertex adjacent in $G$ to an element of $S\})\backslash S$.
If $S=\{v\}$, we will write $g\backslash v$ instead of $g\backslash \{v\}$.
The graph game based on $g$ will be denoted by $\Gamma(g)$ when formality is required, but usually no confusion will arise if we write $g$ instead of $\Gamma(g)$.
The moves in $\Gamma(g)$ are made by selecting vertices of $G$.
The initial move may be of two types:

(i)\ Any available vertex may be selected.

(ii)\ If $G$ has a connected component $C$ which contains no available vertices, then any vertex of $C$ may be selected.

If Player One's initial move is $v$, then $g$ is replaced by $g\backslash v$ with Player Two to move, and the players alternate moves until all vertices have been selected.
The outcome of the game for Player One is 
$$
\sum \{\wt(v)\colon v \text{ was selected by Player One}\}-
\sum \{ \wt(w)\colon w\text{ was selected by Player Two}\},
$$
and the outcome for Player Two is the negative of this.
Of course it is possible to tabulate a score on a running basis during play.
The {\it value} of the game, denoted $\val(g)$, is the outcome for Player One under optimal play.  Note again that this term is not used the same way as in combinatorial game theory.
Of course there is a recursive formula:
$$
\val (g)=\max (\wt(v)-\val (g\backslash v)),\tag1.1
$$
where $v$ ranges through all legal initial moves for $g$ and $G$ is assumed non--empty.

A component $C$ of $G$ is called {\it unbroached} if $C\cap A=\emptyset$ and {\it broached} otherwise, and $g$ is called unbroached if $A=\emptyset$.
Obviously new unbroached components cannot arise in the course of play.
It can be shown that every graph game board $g$ can arise from an unbroached $\tilde g=(\tilde G,\emptyset)$ in optimal play.
Moreover, if $g$ has no unbroached components, $\tilde G$ may be taken connected.
Also if $g$ is a 0--1 game (all weights are 0 or 1) with no unbroached components, then $\tilde g$ may be taken as a 0-1 game with $\tilde G$ connected.

Although all of our theorems, lemmas, etc.~are stated only for graph games as defined above, we will occasionally mention other sorts of games in examples or remarks.
All of these games will be perfect information two--person zero--sum games with real move--by--move scoring and no chance moves, in which the available moves do not depend upon the identity of the player.
For brevity such games will be referred to simply as ``games'' or (for emphasis) ``general games''.
Thus a general game $\Gamma=\Gamma(D)$ is given by an acyclic digraph $D$ with a real weight assigned to each edge, a designated initial node $v_0$, and the requirement that no infinite directed paths starting at $v_0$ exist.
To avoid excess baggage, one should also assume that $v_0$ is the unique source of $D$ and (more) that every node of $D$ can be reached by a directed path starting at $v_0$.
In the play of $\Gamma$ the two players alternately select edges of $D$ so as to construct a directed path starting at $v_0$, and play ends when no more moves are possible.
Each player's score is initially set at 0, and after each move the weight of the selected edge is added to the score of the player who moved and subtracted from the score of the opponent.
The outcome for both players is their final score.
For the graph game $\Gamma(g)$ the nodes of $D$ are the graph game boards $g\backslash S$ which can arise from $g$ in legal play; and for each vertex $v$ which is a legal initial move in $g\backslash S$, there is an edge whose weight is $\wt(v)$ from $g\backslash S$ to $g\backslash (S\cup \{v\})$.
We will almost always assume that $D$ has finitely many nodes.
The only exceptions are the games $\gamma(x)$
 and $\sigma(x)$ in Example 3.20 and the games $\pi(x_1,\dots,x_{\sigma})$ in subsection 4.7.

If $\Gamma_1=\Gamma(D_1)$ and $\Gamma_2=\Gamma (D_2)$ are general games their {\it concatenation}, $\Gamma_1\vee \Gamma_2=\Gamma(D)$, is obtained by placing $\Gamma_1$ and $\Gamma_2$ ``side--by--side'' much like a sum in combinatorial game theory.
Formally, the set of nodes of $D$ is $N_1\times N_2$, where $N_i$ is the set of nodes of $D_i$; for each edge from $v_1$ to $w_1$ in $D_1$ there is an edge of the same weight in $D$ from $(v_1,u_2)$ to $(w_1,u_2)$ for each $u_2$ in $N_2$; and similarly edges in $D_2$ yield edges in $D$.
If $g_1=(G_1,A_1)$ and $g_2=(G_2,A_2)$ are graph game boards, then $g_1\oplus g_2$ is $(G,A)$ where $G$ is the disjoint union of $G_1$ and $G_2$ and $A=A_1\cup A_2$.
It is obvious that the graph game for $g_1\oplus g_2$ is the concatenation of the graph games for $g_1$ and $g_2$; i.e., $\Gamma (g_1\oplus g_2)=\Gamma(g_1)\vee \Gamma(g_2)$ (where we are identifying games which are isomorphic in the obvious sense).

If $G$ has $n$ vertices, then every play of $g=(G,A)$ has exactly $n$ moves.
We will call $n$ the {\it size} of $g$, denoted $\size (g)$, and say that $g$ is {\it even} or {\it odd} according as $n$ is even or odd.
Thus graph games have predetermined parity:\ \ If $g$ is even, then Player Two always makes the last move in $\Gamma(g)$; and if $g$ is odd then Player One always makes the last move.
As discussed more fully in the first paragraph of \S 3,  many of the results in that section apply to all general games with predetermined parity.  In [4] W. Johnson has introduced a class of games, called well-tempered scoring games, in which the available moves are allowed to depend upon the identity of the player.  All general games with predetermined parity can be made into well-tempered scoring games.  Although there is almost no overlap between our \S 3 and [4], it seems worthwhile to investigate the interplay between the ideas of [4] and \S 3, and some such investigation has already taken place.

It may seem unnatural to allow negative weights, but in fact this does no harm and is actually helpful.
To see that it does not harm, denote by $\tau_\lambda(g)$ the result of adding $\lambda$ to each weight.
Clearly for any play of $g$ with outcome $s$, the same moves in $\tau_\lambda(g)$ produce an outcome of $s$ or $s+\lambda$, according as $g$ is even or odd.
Thus optimal strategies are the same for $g$ and $\tau_\lambda(g)$, and $\tau_\lambda(g)$ has positive weights if $\lambda$ is large enough.

Also
$$
\val (\tau_\lambda (g))=\cases \val(g),&\text{$g$ even}\\
\lambda+\val(g),&\text{$g$ odd}.\endcases\tag1.2
$$
In the next section we will solve two--ended stack and multiple stacks games by means of a reduction process.
This process can produce negative weights even when all the original weights are positive.

Some additional notations:\ \ If $g=(G,A)$, $|g|$ denotes $\sum \{ |\wt(v)|\colon v$ is a vertex of $G\}$.
If $G$ is a path whose vertices are $v_1,\ldots,v_n$, in order, and if $\wt(v_i)=a_i$, then $\path (a_1,\ldots,a_n)$ denotes the graph game board $(G,\emptyset)$, $\st(a_1,\ldots,a_n)$ denotes $(G,\{v_1\})$, and $\tes (a_1,\ldots,a_n)$ denotes $(G,\{v_1,v_n\})$.
If $G$ is a cycle whose vertices are $v_1,\ldots,v_n$, in order, and if $\wt(v_i)=a_i$, then $\cyc(a_1,\ldots,a_n)$ denotes $(G,\emptyset)$.
A game of the form $\st (a_1,\ldots,a_n)$ is called a {\it stack}, and the multiple stack games of the introduction are concatenations of stacks.  
The symbol $\langle a_1,\ldots,a_n\rangle$ denotes a graph game board with $n$ vertices, all available, with respective weights $a_1,\ldots,a_n$.
Also $\langle \ \ \rangle$ denotes the empty graph game.

Certain features of $(G,A)$ have no bearing on the graph game and may be altered at will:

1.\ There is no loss of generality in assuming $G$ is a simple graph, since the presence of loops or multiple edges between two vertices is irrelevant.
If a construction done in \S2, say, would produce a non--simple graph $\tilde G$, we assume without mention that $\tilde G$ is replaced with its associated simple graph.

2.\ If a component of $G$ consists of a single vertex $v$, it is irrelevant whether or not $v$ is in $A$.

3.\ If $v$ and $w$ are distinct elements of $A$, it is irrelevant whether or not there is an edge between $v$ and $w$.

\definition{Definition 1.1}We say that $g_1$ is {\it equivalent} to $g_2$, $g_1\sim g_2$, if $\val(g_1\oplus h)=\val (g_2\oplus h)$ for all graph games $h$.
The equivalence class of $g$ will be denoted $[g]$.
Also $[a_1,\ldots,a_n]$ denotes $[\langle a_1,\ldots,a_n\rangle]$ and [\; ] or $\bold 0$ is the equivalence class of the empty game.
Since equivalence is compatible with concatenation, we may define $[g_1]+[g_2]=[g_1\oplus g_2]$.
Then $\sS$, the set of equivalence classes of graph games, is a commutative semigroup, and [\; ] is the identity element of $\sS$.
\enddefinition

Consider the following statements about a graph game $g$ and a real number $x$:
$$
\gathered
\text{(i) }\val(g\oplus h)=\val(h)-x\text{ for all even graph games }h\\
\text{(ii) }\val(g\oplus h)=\val(h)+x\text{ for all odd graph games }h.
\endgathered\tag1.3
$$
It can be shown that (i) and (ii) are equivalent and that (1.3) implies that $g$ is even and $x\geq 0$.
(One way to justify this assertion is to concatenate with $\langle y\rangle$ where either $y>>0$ or $y<<0$.)
If $x\geq 0$, then $\st(a,a+x)$ satisfies (1.3).
This can easily be proved by induction on $\size(h)$, and a generalization will be proved in Proposition 2.1 below.
For $x\geq 0$, the set of all $g$ satisfying (1.3) is an equivalence class which will be denoted by $\ev(x)$.
Of course $\ev(0)=[\; ]$.
The symbol $\ev(x)$ may also denote an unspecified game in the class $\ev(x)$ when no confusion will arise.
(The symbol $\ev$ stands for ``even vertex''.)

\example{Remark 1.2}  Although neither (1.3)(i) nor (1.3)(ii) can hold for all $h$ unless $x\geq 0$ and $g$ is in $\ev(x)$, if $\delta > 0$ and $\{h_i\}$ is a finite set of graph games, then there is a game $g=\langle x,x+\delta\rangle$ such that $\val(g\oplus h_i) = \val(h_i)+\delta$, $\forall i$.  There is also a game $g'= \langle y,y+\delta\rangle$ such that $\val(g'\oplus h_i)= \val(h_i) + \delta$ for $h_i$ even, and $\val(g'\oplus h_i)=\val(h_i) -\delta$ for $h_i$ odd.
  (Just take $x>>0$ and $y<<0$.)  It can then be shown that if $g_1$ and $g_2$ are inequivalent graph games, there is a graph game $h$ such that $\val(g_1\oplus h)$ and $\val(g_2\oplus h)$ have opposite signs.  Moreover, one can specify the parity of $h$ or one can specify which of $\val(g_i\oplus h)$ is positive.  However, one cannot in general specify both simultaneously.  (But if either $g_1$ and $g_2$ are both invertible or if $g_1$ and $g_2$ have opposite parities, one can specify both simultaneously.)
\endexample

For graph games $g_1,g_2$ let $\sI([g_1],[g_2])$ or $\sI(g_1,g_2)$ denote $\sup|\val (g_1\oplus h)-\val (g_2\oplus h)|$, where the supremum is extended over all graph games $h$.
Then $\sI([g_1],[g_2]$ may be $+\infty$, but it satisfies all other requirements to be a metric on $\sS$.
In particular $\sI(g_1,g_2)=0$ if and only if $g_1\sim g_2$.
The following elementary lemma is stated for the purpose of reference.

\proclaim{Lemma 1.3}If $g'$ is obtained from $g$ by changing one weight from $a$ to $a'$, then $\sI(g,g')\leq |a-a'|$.
\endproclaim

\proclaim{Proposition 1.4}For graph games $g_1$ and $g_2$, $\sI(g_1,g_2)<\infty$ if and only if $g_1$ and $g_2$ have the same parity.
\endproclaim

\demo{Proof}It is easy to prove by induction that if all weights in $g$ are 0, then $g$ is equivalent to $\langle\ \ \rangle$ or $\langle 0\rangle$ according as $g$ is even or odd.
It then follows from Lemma 1.3 that $\sI(g_1,g_2)\leq |g_1| + |g_2|$ if $g_1$ and $g_2$ have the same parity.
If $g_1$ and $g_2$ have opposite parity, then $|\val (g_1\oplus h)-\val (g_2\oplus h)|$ can be made arbitrarily large by taking $h$ in $\ev(x)$ for $x>>0$ or $h=\langle y\rangle$ for $y<<0$.
\enddemo

If $g$ and $h$ have the same parity, $v$ is a legal initial move in $g$ and $w$ is a legal initial move in $h$, we say that the $v$--move {\it metrically dominates} the $w$--move if $\sI(g\backslash v, h\backslash w)\leq \wt(v)-\wt(w)$.
It is allowed that $g=h$, and it is clear that metric domination is a transitive relation.
If each move metrically dominates the other, we say that they are {\it equivalent}.
The ``metric domination principle'' is that if the $v$--move in $g$ metrically dominates the $w$--move in $h$, then the $v$--move in $g\oplus k$ will produce at least as large an outcome for Player One as the $w$--move in $h\oplus k$ (assuming optimal play after the initial moves).
(The converse is also true.)
For example, it follows from Lemma 1.3 that if $a_1\geq a_2,\ldots,a_n$, then the $a_1$--move metrically dominates all others in $\langle a_1,\ldots,a_n\rangle$.
This implies the (obvious) fact that the greedy strategy is optimal in this game.

The following facts are easily proved:
$$
\text{If }a_1\geq a_2\geq \ldots a_n, \text{ then }\val (\langle a_1,\ldots,a_n\rangle)
=\sum_1^n (-1)^{i-1} a_i.  \tag1.4
$$
$$
 [a,a]=[\; ]\text{ and }[a, a, a_1,\ldots,a_n]=[a_1,\ldots,a_n].\tag1.5
$$
$$
 \ev(x)+\ev(y)=\ev(x+y).\tag1.6
$$

Let $\sG_0=\{ [a_1,\ldots,a_n]\colon n\geq 0, a_1,\ldots,a_n\in\bR\}$, where $\bR$ is the set of real numbers, and let $\sS_0=\{[g]+\ev(x)\colon [g]\in\sG_0$ and $x\geq 0\}$.
Then $\sG_0$ is a group (by (1.5)), and $\sS_0$ is a semigroup.
Although these objects seem trivial, it will be shown in Corollary 2.12 below that $\sS_0$ can also be described either as the set of all equivalence classes of stacks or as the set of all equivalence classes of trees with one available vertex, and $\sG_0$ can also be described as the set of all equivalence classes of two--ended stacks.
In \S3 the completion of $\sG_0$, with respect to $\sI$, will be described in terms of Lebesgue measure on $\bR$.

\proclaim{Lemma 1.5 (limited badness principle)}Let $g=(G,A)$ be a graph game with no unbroached components.
Let $u,w\in A,u\neq w$, $\wt(u)=a$, and $\wt(w)=b$.
Let $m=\max \{\wt (v):v$ any vertex$\}$ and $m'=\max\{\wt(v)\colon v$ is a vertex and $v\neq u\}$.
Then:
\itemitem{(i)}We have $\val(g)\leq a-\val(g\backslash u)+2(m-a)$
\itemitem{(ii)}We have $b-\val(g\backslash w)\leq a-\val (g\backslash u)+2(m'-a)$.
\endproclaim

\demo{Proof}Use induction on $n=\size(g)$, the case $n=2$ being obvious.
Clearly (i) for $n$ follows from (ii) for $n$, and we will deduce (ii) for $n$ from (i) for $n-1$.
Note that $\val (g\backslash w)\geq a-\val (g\backslash \{u,w\})$, and therefore the left side of (ii) is at most $b-a+\val(g\backslash \{u,w\})$.
By induction $\val(g\backslash u)\leq b-\val (g\backslash \{u,w\})+2(m'-b)$, and therefore the right side of (ii) is at least $a-b+\val (g\backslash \{u,w\})-2 (m'-b) +2(m'-a)$.
The result follows.
\enddemo

\proclaim{Corollary 1.6}Assume the graph game $g$ has no unbroached components and $v$ is an available vertex.
If $\wt(v)=\max \{\wt (u)\colon u\text{ any vertex}\}$, then the $v$--move is optimal.
If also, $\wt(u) < \wt(v)$ for every vertex $u$ different from $v$, then the $v$--move is the only optimal initial move.
\endproclaim

\example{Example 1.7}We show that the hypothesis of no unbroached components cannot be eliminated in Corollary 1.6 or Lemma 1.5.
If $g_1=\langle m\rangle +\cyc (0,1,0,1,\ldots,1,0,1,1)$, where $m\geq 1$, the cycle contains $k$ 1--vertices, and $k>2m$, then the only optimal moves are the two ``rightmost'' 1--vertices in the cycle (the value of $g_1$ is $k-m$).
If $k>>m$, the cost of initially making the $m$--move is very large.
If $g_2=\langle 0\rangle \oplus g_1$ and $k> 2m+2$, the only optimal move is one of the fairly central 1--vertices in the cycle, and the value is $-m$ or $1-m$ according as $k$ is even or odd.
Again the cost of the $m$--move can be very high.
The analysis of $g_1$ and $g_2$ is most easily done with the help of the theory presented later in the paper.
The reason for choosing these examples is that $g_1$ and $g_2$ are both invertible.
Thus invertibility combined with a parity hypothesis cannot rescue the principle of limited badness if there is an unbroached component.
Easier examples with non--invertible games are $h_1=\langle m\rangle \oplus \path (1,0,1,0,\ldots,1,0)$ where there are $k$ 1--vertices in the path and $k> m\geq 1$ and $h_2=\langle 0\rangle \oplus h_1$.
(Both $h_1$ and $h_2$ have value $k-m$.)
Also note that for the games $\path (1,0,0)$ (which is equivalent to $\langle 1\rangle$) and $\cyc(0,1,0,0)$ (which is equivalent to $\langle 1,0\rangle$), the 1--move is optimal but it is not the only optimal initial move.

There are also examples to show that the second assertion of Corollary 1.6 would fail if we assume only that $\wt(v)$ is maximal among all vertices and $\wt(v) > \wt(u)$ for all available vertices $u$ other than $v$.  In each of the games $\langle 0,1\rangle\oplus\st(0,1,0,-1)$, $\langle 0,1\rangle\oplus\st(0,1,0,-1,0)$, $\tes(3,2,2,4,4)$, and $\tes(3,2,1,2,4,4)$ there are two optimal initial moves.
\endexample

\bigskip\noindent
{\bf \S2.\ Stacks and two--ended stacks}.

The three main reductions used in solving concatenations of stacks and two--ended stacks can be applied to somewhat more general graph games.
Therefore they are given in abstract form.

The following condition on a finite sequence, $(a_1,\ldots,a_n)$ of real numbers will often be used in this section.
$$
\sum^{2k}_{i=1} (-1)^{i-1} a_i\leq 0\text{ for }1\leq k\leq \lfloor{n\over 2}\rfloor.\tag"(ASC)"
$$
In particular $\ASC$ is vacuously true if $n=1$.
We say that $(a_1,\ldots,a_n)$ is an {\it $\ev$--sequence} if $n$ is even and $\ASC$ holds.
In this case the {\it weight} of the \ev--sequence is the number $x=\sum\limits^n_{i=1} (-1)^i a_i$.
Note that $x\geq 0$.

\proclaim{Proposition 2.1}If $(a_1,\ldots,a_n)$ is an $\ev$--sequence of weight $x$, then $\st(a_1,\ldots,a_n)$ is in $\ev(x)$.
\endproclaim

\demo{Proof}We use induction on $n$.
The case $n=2$ was claimed as easy in \S1, but its proof will be included below.
If $g=\st (a_1,\ldots,a_n)$ and $h$ is any graph game, let
$$
v^*=\cases\val(h)-x,\text{if $h$ is even}\\
\val(h)+x,\text{if $h$ is odd} \endcases .
$$
We show by induction on the size of $h$ that $\val(g\oplus h)=v^*$.
The case where $h$ is empty is obvious.
If $h$ is non--empty, then it follows from the induction on $\size(h)$ that
$\sup(\wt (v)-\val (g\oplus h\backslash v))=v^*$, where the supremum is taken over all legal initial moves $v$ in $h$.
Thus it is sufficient to show that $a_1-\val (\st (a_2,\ldots,a_n) \oplus h)\leq v^*$.
But clearly $a_1-\val (\st (a_2,\ldots,a_n)\oplus h)\leq a_1 -a_2 +\val (\st (a_3,\ldots,a_n)\oplus h)$.
If $n=2$, we are done.
If $n\geq 4$, note that $(a_1-a_2+a_3,a_4,\ldots,a_n)$ is again an $\ev$--sequence of weight $x$.
By the induction on $n$, $g'$ is in $\ev(x)$ for $g'=\st(a_1-a_2+a_3,a_4,\ldots,a_n)$.
Also, by Lemma 1.3, $\val (\st(a_3,\ldots,a_n)\oplus h)\leq a_2-a_1+\val (g'\oplus h)=a_2-a_1+v^*$.
\enddemo

\proclaim{Theorem 2.2 (detaching $\ev$--sequences)}Assume $(a_1,\ldots,a_n)$ is an $\ev$--sequence of weight $x$.
Let $g=(G,A)$ be a graph game board such that $G$ contains the graph of $\path (a_1,\ldots,a_n)$ as a complete subgraph with corresponding vertices denoted by $v_1,\ldots,v_n$.
Assume there is an edge joining $v_1$ to a vertex $w$ in $G\backslash \{v_1,\ldots,v_n\}$, and there are no other edges between $\{v_1,\ldots,v_n\}$ and $G\backslash \{v_1,\ldots,v_n\}$.
Finally, assume that $A\cap \{v_1,\ldots,v_n\}=\emptyset$ and the component of $G$ containing $w$ is broached.
Let $\tilde g$ be the graph game board obtained from $g$ by removing $\{v_1,\ldots,v_n\}$ and the incident edges without changing $A$.
Then $[g]=[\tilde g]+\ev(x)$.
\endproclaim

\demo{Proof}For any graph game $h$, let $v^*=\cases \val(\tilde g\oplus h)-x,\text{ if $\tilde g\oplus h$ is even}\\ \val(\tilde g\oplus h)+x,\text{if $\tilde g\oplus h$ is odd}\endcases$.
We show by induction on $\size(\tilde g\oplus h)$ that $\val(g\oplus h)=v^*$.
The case where $\size(\tilde g\oplus h)=1$ is obvious, since then $w$ is the only vertex in $\tilde g\oplus h$.
Note that there is a one--one correspondence between the legal initial moves in $g\oplus h$ and $\tilde g\oplus h$.
For any such move $u$, if $u\neq w$ the induction hypothesis applies to $(g\oplus h)\backslash u$.
And if $w\in A$, Prop.~2.1 applies to $(g\oplus h)\backslash w$, since $g\backslash w=(\tilde g\backslash w)\oplus \st (a_1,\ldots,a_n)$.
\enddemo

\proclaim{Corollary 2.3}Under the same hypotheses $g\sim \tilde g\oplus \st (a_1,\ldots,a_n)$.
\endproclaim

\proclaim{Corollary 2.4}If $(a_1,\ldots,a_n)$ is an $\ev$--sequence, then
\newline $\st(b_1,\ldots,b_m,a_1,\ldots,a_n)\sim \st (b_1,\ldots,b_m)\oplus \st (a_1,\ldots,a_n)$.
\endproclaim

\proclaim{Theorem 2.5 (detaching available maxima)}If $v$ is an available vertex in the graph game $g$, and if $\wt(v)$ is maximal among all the vertices of $g$, then $g\sim (g\backslash v)\oplus \langle \wt(v)\rangle$.
Moreover, the $v$--move in $g$ metrically dominates all other moves of available vertices in $g$ or in $(g\backslash v)\oplus \langle \wt(v)\rangle$.
\endproclaim

\demo{Proof}We use induction on $\size(g)$, the case $\size(g)=1$ being trivial.
Let $m=\wt(v)$.
We first prove the metric domination claim.
Clearly the $v$--move in $g$ is equivalent to the $\langle m\rangle$--move in $(g\backslash v)\oplus \langle m\rangle$.
If $u$ is any other available vertex in either game, then $u$ is available in $g\backslash v$.
Let $g'$ be the game obtained from $g\backslash v$ by changing the weight of $u$ to $m$.
Then by induction $g'\sim (g\backslash \{ u,v\})\oplus \langle m\rangle$, and if $u$ is available in $g$ also, then $g\backslash u\sim (g\backslash \{u,v\})\oplus \langle m\rangle$.
Also by Lemma 1.3 $\sI(g',g\backslash v)\leq m-\wt(u)$.
These facts imply the metric domination.
\enddemo

Next we prove that $\val(g\oplus h)=\val ((g\backslash v)\oplus \langle m\rangle \oplus h)$, for all graph games $h$.
If this is false choose a counterexample with $\size(h)$ as small as possible.
Then let a vertex $u$ be an optimal move in the game with the larger value.
It is impossible for $u$ to be in $h$ because then the minimality assumption on $\size(h)$ would imply that the $u$--move leads to the same outcome in both games.
Similarly, $u$ cannot be a legal move in $g$ other than $v$, because any such move is legal in both games and the induction on $\size(g)$ would apply.
Obviously $u$ cannot be $v$ or the vertex of $\langle m\rangle$.
The remaining possibility is that $u$ is a vertex in $g\backslash v$, adjacent
 to $v$ in $g$.
But this case is ruled out by the already proved metric domination.

\example{Remarks}The context of this theorem is somewhat similar to that of Corollary 1.6, but neither result implies the other.
In Theorem 2.5, $g$ may have unbroached components, but they play no real role.
(They may as well be consolidated with the game $h$ that appears in the proof, and of course $h$ may have vertices of weight larger than $m$.)
\endexample

\proclaim{Theorem 2.6 (condensation)}Assume that the graph game $g$ has three distinct vertices $u,v,w$ such that:

1.\ The weight of $v$ is maximal among all vertices of $g,v$ is not available, and $v$ is adjacent to $u$ and $w$ and no other vertices;

2.\ Each of $u,w$ has degree at most 2 and has degree 1 if available;

3.\ The component containing $u,v,w$ is broached.

Let $\tilde g$ be the game obtained by condensing $u,v$, and $w$ to a single vertex $p$ as follows:

$1'$.\ A vertex $z$ other than $u,v,w$ is adjacent to $p$ in $\tilde g$ if and only if it is adjacent to $u$ or $w$ in $g$;

$2'$.\ The vertex $p$ is available in $\tilde g$ if and only if $u$ or $w$ is available in $g$;

$3'$.\ The weight of $p$ is $\wt(u)-\wt(v)+\wt(w)$.

Then $\tilde g$ is equivalent to $g$.
Moreover if $u$ (or $w$) is available in $g$, then the $u$--move (or $w$--move) metrically dominates the $p$--move in $\tilde g$.
\endproclaim

\demo{Proof}Let $a=\wt(u),\ m=\wt(v),\ b=\wt(w)$, and $c=a-m+b=\wt(p)$.

We first prove the metric domination.
If, say, $u$ is available in $g$, then by Theorem 2.5, $g\backslash u\sim (g\backslash \{u,v\})\oplus \langle m\rangle$.
Let $g'$ be the game obtained from $g\backslash \{u,v\}$ by changing the weight of $w$ to $m$.
Then by Lemma 1.3, $\sI (g\backslash \{u,v\}, g')\leq m-b$, and by Theorem 2.5 $g'\sim (g\backslash \{u,v,w\})\oplus \langle m\rangle=(\tilde g\backslash p)\oplus \langle m\rangle$.
From this and formula (1.6) we see that $\sI (g\backslash u,\tilde g\backslash p)\leq m-b=a-c=\wt(u)-\wt(p)$.

Next we prove that if $u$ (or $w$) is available, then for any game $h$, $a-\val ((g\backslash u)\oplus h)\leq \val(\tilde g\oplus h)$ 
(or $b-\val ((g\backslash w)\oplus h)\leq \val(\tilde g\oplus h))$.
Clearly $a-\val ((g\backslash u)\oplus h)\leq a-m+\val ((g\backslash \{u,v\})\oplus h)$.
But by Lemma 1.3 $\sI(g\backslash \{u,v\},\tilde g)\leq |b-c|=m-a$.

Now to prove the equivalence we need to show that $\val (g\oplus h)=\val (\tilde g\oplus h)$ for all $h$.
If this is false choose a counterexample such that $\size(g)+\size(h)$ is as small as possible.
Let $z$ be a vertex which is an optimal move in the game with the larger value.
Then by the metric domination principle $z$ cannot be $p$, and by the paragraph above $z$ cannot be $u$ or $w$.
But this implies that $z$ is a legal move in both games, and the minimality hypothesis implies that $\val ((g\oplus h)\backslash z)=\val ((\tilde g\oplus h)\backslash z)$, a contradiction.
\enddemo

\example{Remark}If $z$ is a legal move other than $p$ in $\tilde g$, then $z$ is also a legal move in $g$; and Theorem 2.6 implies that the $z$--moves in $g$ and $\tilde g$ are equivalent.
Thus each legal move in $\tilde g$ is metrically dominated by a corresponding move in $g$.
\endexample

\proclaim{Corollary 2.7}Let $g$ be $\st (a_1,\ldots,a_n)$ or $\tes$ $(a_1,\ldots,a_n)$ and assume $a_i=\max (a_1,\ldots,a_n)$ where $1< i< n$.
If $\tilde g$ is respectively $\st(a_1,\ldots,a_{i-2},a_{i-1}-a_i+a_{i+1},a_{i+2},\ldots,a_n)$ or $\tes(a_1,\ldots,a_{i-2},a_{i-1}-a_i+a_{i+1},a_{i+2},\ldots,a_n)$, then $\tilde g\sim g$.
Moreover each legal move in $g$ (left or right) metrically dominates the corresponding move in $\tilde g$.
\endproclaim

The above results can be used to reduce any concatenation of stacks and two--ended stacks to an equivalent game for which an optimal strategy is obvious.
Moreover, one could deduce an optimal initial move for the original game from one for the reduced game.
If we made a choice of the rule to be used when more than one reduction is possible, we would have an algorithm for solving such games; but it would be a somewhat slow algorithm ($O(n^{2})$ steps if $n=\size(g)$).
A fast ($O(n)$ steps) algorithm will be presented below.
The easiest way to explain and justify this algorithm is to leap directly to a description of the ``fully condensed'' game.

A finite sequence $(a_1,\ldots,a_n)$ is called a {\it slice} if $n$ is odd and both $(a_1,\ldots,a_n)$ and its reverse, $(a_n,\ldots,a_1)$, satisfy $\ASC$.
If $S=(a_1,\ldots,a_n)$ is a slice its {\it weight} or {\it slice--weight}, denoted by $\wt(S)$, is the number $s=\sum_1^n (-1)^{i-1} a_i$.
Note that $s\leq a_1$ and $s\leq a_n$.

\proclaim{Lemma 2.8}If the finite sequence $(a_1,\ldots,a_n)$ is partitioned into slices $S_1,\ldots,S_m$, where $s_i=\wt(S_i)$, then $(a_1,\ldots,a_n)$ satisfies $\ASC$ if and only if $(s_1,\ldots,s_m)$ satisfies $\ASC$.
It follows that $(a_1,\ldots,a_n)$ is a slice or $\ev$--sequence if and only if $(s_1,\ldots,s_m)$ is.
In this case both sequences have the same weight.
\endproclaim

\demo{Proof}1.\ It is obvious that $\ASC$ for $(a_1,\ldots,a_n)$ implies $\ASC$ for $(s_1,\ldots,s_m)$, since each of the alternating sums appearing in $\ASC$ for $(s_1,\ldots,s_m)$ is equal to one of the alternating sums in $\ASC$ for $(a_1,\ldots,a_n)$.

2.\ Now assume $\ASC$ for $(s_1,\ldots,s_m)$ and let $1\leq k\leq \lfloor{n\over 2}\rfloor$.
First suppose that $a_{2k}$ falls in $S_\ell$ where $\ell$ is odd.
Let $a_p$ be the first term in $S_\ell$ ($p$ is odd).
Then $\sum\limits^{2k}_{i=1} (-1)^{i-1} a_i=\sum\limits^{\ell-1}_{j=1} (-1)^{j-1} s_j+\sum\limits^{2k}_{i=p} (-1)^{i-1} a_i\leq 0+0=0$, where the inequalities follow from $\ASC$ for $(s_1,\ldots,s_m)$ and $S_\ell$.
Next suppose that $a_{2k}$ falls in $S_\ell$ for $\ell$ even.
Let $a_q$ be the last term in $S_\ell$ ($q$ is even).
Then $\sum\limits^{2k}_{i=1} (-1)^{i-1} a_i=\sum\limits^\ell_{j=1} (-1)^{j-1} s_j+\sum\limits^{q-2k}_{r=1} (-1)^{r-1} a_{q+1-r}\leq 0+0=0$, where the inequalities follow from $\ASC$ for $(s_1,\ldots,s_m)$ and the reverse of $S_\ell$.

The rest of the lemma is obvious.
\enddemo

\definition{Definition 2.9}A sequence $(s_1,\ldots,s_m)$ will be called {\it $U$--shaped} if it satisfies the following equivalent conditions:

(i)\ If $1< j< m$, then either $s_j < s_{j+1}$ or $s_j < s_{j-1}$ (no interior (weak) relative maxima).

(ii)\ Let $k$ be the smallest index such that $s_k=\min (s_1,\ldots,s_m)$ and $\ell$ the largest index such that $s_\ell=\min (s_1,\ldots,s_m)$.
Then $\ell=k$ or $\ell=k+1$, $(s_1,\ldots,s_k)$ is strictly decreasing, and $(s_\ell,\ldots,s_m)$ is strictly increasing.

Also, $(s_1,\ldots,s_m)$ will be called {\it weakly $U$--shaped} if either $m=1$ or there is $p$ such that $1\leq p\leq m-1$, $(s_1,\ldots,s_p)$ is monotone non--increasing, and $(s_{p+1},\ldots,s_m)$ is monotone non--decreasing.

Note that $U$--shaped implies weakly $U$--shaped, $m\leq 2$ implies $(s_1,\ldots,s_m)$ is $U$--shaped, and $(s_1,\ldots,s_m)$ $U$--shaped or weakly $U$--shaped implies the same for $(s_2,\ldots,s_m)$ and $(s_1,\ldots,s_{m-1})$ if $m\geq 2$.
\enddefinition

Part (i) of the following lemma is intended for two--ended stacks and part (ii) for stacks.

\proclaim{Lemma 2.10}Let $(a_1,\ldots,a_n)$ be a finite sequence of real numbers.
Then 

(i) The sequence $(a_1,\ldots,a_n)$ can be partitioned into slices, $S_1,\ldots,S_m$ such that the sequence $(s_1,\ldots,s_m)$ of slice--weights is $U$--shaped.

(ii)\ The sequence $(a_1,\ldots,a_n)$ can be partitioned into $S_1 S_2\ldots S_m X$, where possibly $X$ is missing or $m=0$, such that $S_1,\ldots,S_m$ are slices, $X$ is an $\ev$--sequence, and the sequence $(s_1,\ldots,s_m)$ of slice--weights is strictly decreasing.
\endproclaim

\demo{Proof}(i)\ We use induction on $n$, the cases $n=1$ or 2 being trivial.
If $(a_1,\ldots,a_n)$ is not already $U$--shaped, choose $i$ such that $1< i<n$ and $a_i\geq a_{i-1},a_{i+1}$.
Then the induction hypothesis implies a suitable partition of $(a_1,\ldots,a_{i-2},a_{i-1}-a_i+a_{i+1},a_{i+2},\ldots,a_n)$ into slices.
Since $(a_{i-1},a_i,a_{i+1})$ is a slice, Lemma 2.8 implies that $(a_1,\ldots,a_n)$ has a corresponding partition into slices.

(ii)\ Let $S'_1,\ldots,S'_p$ be the slices found in part (i), and let $k$ be as in condition (ii) of Definition 2.9.
If $k=p$, we are done.
If $k<p$ and $p-k$ is odd, then let $m=k-1$, $S_j=S'_j$ for $j\leq m$ (if $m>0$), and $X=S'_k S'_{k+1}\ldots S'_p$.
If $k<p$ and $p-k$ is even, then let $m=k$, $S_j=S'_j$ for $j\leq m$, and $X=S'_{k+1}\ldots S'_p$.
In either case Lemma 2.8 implies that $X$ is an $\ev$--sequence.
\enddemo

\example{Remark}In part (ii) let $x$ be the weight of $X$ if $X$ is present and let $x=0$ if $X$ is missing.
(Of course $x$ might be 0 even if $X$ is present.)
Then it turns out that the slice--weight sequence $(s_1,\ldots,s_m)$ in part (i)  and the object $(s_1,\ldots,s_m;x)$ in part (ii) are unique.
This will follow from the justification of the algorithm, 2.18 below.
However, the actual partitions into slices may not be unique.
\endexample

\proclaim{Theorem 2.11}Let $(a_1,\ldots,a_n)$ be a finite sequence of real numbers.

(i)\ If $(a_1,\ldots,a_n)$ is partitioned into $S_1 S_2\ldots S_m$, where the $S_i$'s are slices and the sequence $(s_1,\ldots,s_m)$ of slice--weights is weakly $U$--shaped, then $\tes (a_1,\ldots,a_n)\sim \langle s_1,\ldots,s_m\rangle$.

(ii)\ If $(a_1,\ldots,a_n)$ is partitioned into $S_1 S_2 \ldots S_m X$, where possibly $m=0$ or $X$ is missing, such that $S_i$ is a slice of weight $s_i$, $X$ is an $\ev$--sequence of weight $x$, and the sequence $(s_1,\ldots,s_m)$ is non--increasing, then $\st(a_1,\ldots,a_n)\sim \langle s_1,\ldots,s_m\rangle \oplus\ev(x)$.
\endproclaim

\demo{Proof}(i)\ If this is false, choose a counterexample $(a_1,\ldots,a_n)$ with $n$ as small as possible.
Let $M=\max (a_1,\ldots,a_n)$.

Case 1. There is a slice $S_i$ of size at least 3 such that $a_j=M$ for an interior point $a_j$ of $S_i$.
In this case Corollary 2.7 implies $\tes(a_1,\ldots,a_n)\sim \tes(a_1,\ldots,a_{j-2},a_{j-1}-a_j+a_{j+1},a_{j+2},\ldots,a_n)$.
Moreover, the new sequence is partitioned into slices with the same weight sequence $(s_1,\ldots,s_m)$.
Then the minimality of $n$ implies the desired equivalence, a contradiction. 

Case 2. Case 1 does not hold and either $a_1=M$ or $a_n=M$.
If $a_1=M$ and $S_1$ has size at least 3, then also $a_2=M$, a contradiction.
So $S_1$ consists only of $a_1$ and $s_1=M$.
By Theorem 2.5 $\tes(a_1,\ldots,a_n)\sim \langle s_1\rangle+\tes(a_2,\ldots,a_m)$.
Also $(a_2,\ldots,a_m)$ is partitioned into $S_2\ldots S_m$ with weight sequence $(s_2,\ldots,s_m)$.
So the minimality of $n$ implies $\tes(a_2,\ldots,a_m)\sim \langle s_2,\ldots,s_m\rangle$, establishing the desired equivalence.

Case 3. Neither Case 1 nor Case 2 holds.
Then $a_j=M$ for some $j$ with $1< j<n$ and $a_j$ is an endpoint of some $S_i$.
If $S_i$ has size at least 3, then the adjacent interior point of $S_i$ is also $M$, a contradiction.
Thus $S_i$ has size 1 and $s_i=M$.
Since $(s_1,\ldots,s_m)$ is weakly $U$--shaped and both endpoints of $S_k$ are at least $s_k$, this implies either $s_1=s_2=\ldots=s_i=M$ or $M=s_i=s_{i+1}=\ldots=s_m$.
This implies $a_1=M$ or $a_n=M$.
So Case 3 cannot occur.

(ii)\ The proof is similar.
Case 1 is modified by also including the possibility that $M$ occurs at an interior point of $X$.
In Case 2, if $a_n=M$, we use Proposition 2.1 and Corollary 2.4 to deduce $\st(a_1,\ldots,a_n)\sim\st (a_1,\ldots,a_{n-2})\oplus\ev (M-a_{n-1})$.
If $(a_{n-1},a_n)$ is part of $X$, we are done, since $X\backslash (a_{n-1},a_{n-2})$ is still an $\ev$--sequence.
Otherwise $a_n$ is the final point of $S_m$ and as above we see that $S_m$ has size 1 and $s_m=M$.
The monotonicity now implies that $s_1=M$, and as above $a_1=M$.
The proof that Case 3 cannot hold is almost identical to the above argument.
We also have to consider the possibility that $M$ occurs as the left endpoint of $X$, and this presents no difficulty.
\enddemo

\proclaim{Corollary 2.12}

(i)  The semigroup $\sS_0$ is the set of all equivalence classes of stacks.

(ii)  The semigroup $\sS_0$ is the set of all equivalence classes of graph game boards $g=(G,A)$, where $G$ is a tree and $A =\{v\}$ for a vertex $v$.

(iii) The group $\sG_0$ is the set of all equivalence classes of two-ended stacks.
\endproclaim

\demo{Proof}It follows from Lemma 2.10 and Theorem 2.11 that $[g] \in\sS_0$ if $g$ is a stack and $[g] \in\sG_0$ if $g$ is a two--ended stack.  Any element $[g]$ of $\sG_0$ may be written as $[s_1,\dots,s_m]$ where the sequence $(s_1,\dots,s_m)$ is (strictly) decreasing, and then $g \sim \tes(s_1,\dots,s_m)$.  Also an element $[g]$ of $\sS_0$ may be written as $[s_1,\dots,s_m] + \ev(x)$ where $(s_1,\dots,s_m)$ is decreasing, and then $g \sim \st(s_1,\dots,s_m,0,x)$.

We prove the remaining part of (ii) by induction on $\size(g)$, the case of size 1 being obvious.  For $v$ the available vertex, $g\backslash v = h_1\oplus\dots\oplus h_k$, where each $h_i$ satisfies the same hypothesis as $g$ and is of smaller size.  Since $\sS_0$ is closed under addition, $[g\setminus v] \in \sS_0$.  Let $g\backslash v \sim \st(a_1,\dots,a_n)$.  Then $g \sim \st(\wt(v),a_1,\dots,a_n)$.
\enddemo

\example{Remark 2.13}  Say that two finite sequences of real numbers are {\it similar} if their slice--weight sequences, as in Lemma 2.10(i) are the same.  What does similarity really mean?  Hypothetical applications of Lemma 2.10 to other areas of mathematics might give a good answer to this question.  We can give a partial answer.  If $\underline s$ is the slice--weight sequence for $(a_1,\ldots,a_n)$, then the equivalence classes $[\tes(a_1,\ldots,a_n)]$, $[\st(a_1,\ldots,a_n)]$, and $[\st(a_n,\ldots,a_1)]$ can all be calculated from $\underline s$.  If $n$ is odd the converse is true:  The slice--weight sequence $\underline s$ can be calculated from the three equivalence classes.  But this is not so if $n$ is even.  If $(b_1,\ldots,b_k)$ is strictly decreasing, $(c_1,\ldots,c_{\ell})$ is strictly increasing, $k$ and $\ell$ are even (possibly 0), and if $x<b_k$ and $x<c_1$, then all of the slice--weight sequences $(b_1,\ldots,b_k,x,x,c_1,\ldots,c_{\ell})$ yield the same three equivalence classes, as does $(b_1,\ldots,b_k,c_1,\ldots,c_{\ell})$.  A possibly related fact is that in the even case, but not in the odd case, the equivalence class of the two--ended stack is determined by the equivalence classes of the two stacks.
\endexample

The next corollary strengthens, the first part of Corollary 2.7 by allowing a much more general type of condensation.

\proclaim{Corollary 2.14}Consider a sequence $(a_1,\ldots,a_n)$ of real numbers and assume $(a_i,\ldots,a_j)$ is a slice of weight $s$.
Then $g_1=\tes (a_1,\ldots,a_n)\sim\tes (a_1,\ldots,a_{i-1},s,a_{j+1},\ldots,a_n)=g'_1$ and $g_2=\st (a_1,\ldots,a_n)\sim\st (a_1,\ldots,a_{i-1},s,a_{j+1},\ldots,a_n)=g'_2$.
\endproclaim

\demo{Proof}Use Lemma 2.10 to partition the sequence $(a_1,\ldots,a_{i-1},s,a_{i+1},\ldots,a_n)$ into $S'_1\ldots S'_m$ or $S''_1\ldots S''_m X$.
By Lemma 2.8 this leads to a corresponding decomposition of $(a_1,\ldots,a_n)$ into $S_1\ldots S_m$ or $S_1\ldots S_m X$ with the same weights.
The desired equivalence follows by applying the theorem to both the original sequence and the condensed sequence.
\enddemo

\proclaim{Corollary 2.15}If $a_1\geq a_2$, then $\tes (a_1,\ldots,a_n)\sim \langle a_1\rangle\oplus\tes (a_2,\ldots,a_n)$ and $\st(a_1,\ldots,a_n)\sim\langle a_1\rangle \oplus\st (a_2,\ldots,a_n)$.
\endproclaim

\demo{Proof}Partition $(a_2,\ldots,a_n)$ as in the hypothesis of part (i) or part (ii) of the theorem.
Since $a_2\geq s_1$ if $m>0$, then the sequence $(a_1,s_1,\ldots,s_m)$ still satisfies the hypothesis of the theorem.
So we can just apply the theorem both to the original sequence and to $(a_2,\ldots,a_m)$.
\enddemo

\proclaim{Corollary 2.16}

(i)\ Under the hypotheses of part (i) of the theorem, the $a_1$--move in $\tes(a_1,\ldots,a_n)$ metrically dominates the $s_1$--move in $\langle s_1,\ldots,s_m\rangle$ and the $a_n$ move in $\tes(a_1,\ldots,a_n)$ metrically dominates the $s_m$--move in $\langle s_1,\ldots,s_m\rangle$.

(ii)\ Under the hypotheses of part (ii) of the theorem, if $m>0$, then the $a_1$--move in $\st(a_1,\ldots,a_n)$ metrically dominates the $s_1$--move in $\langle s_1,\ldots,s_m\rangle\oplus\ev(x)$.
\endproclaim

\demo{Proof}(i)\ It is enough to prove this for the $a_1$--move.
We must show $\sI(\tes(a_2,\ldots,a_n)$, $\langle s_2,\ldots,s_m\rangle)\leq a_1-s_1=\sum_2^j (-1)^i a_i\geq 0$ if $S_1=(a_1,\ldots,a_j)$.
This is trivial if $j=1$.
If $j>1$, we replace $a_2$ by $a'_2=a_2-(a_1-s_1)=\sum^j_3 (-1)^{i-1} a_i$, and show $\tes(a'_2,a_3,\ldots,a_n)\sim \langle s_2,\ldots,s_m\rangle$.
Note that $\ub=(a'_2,a_3,\ldots,a_j)$ is an even length sequence whose reverse satisfies $\ASC$ and whose alternating sum is 0.
It follows that $\ub$ is an $\ev$--sequence of weight 0.
This implies that $\ub$ is of the form $S' S''$ where $S'$ and $S''$ are slices of the same weight $s'$.
To find $S'$ and $S''$, just choose $S'$ to be an initial odd segment whose alternating sum is minimal.
The fact that any initial odd segment of $S'$ has an alternating sum that is at least as large implies that the reverse of $S'$ satisfies $\ASC$.
And the fact that $S'Y$ has an alternating sum that is at least as large, where $Y$ is an initial even segment of $S''$, implies that $S''$ satisfies $\ASC$.
So $S'$ and $S''$ are slices, and it it obvious that they have the same weight.
Moreover $s'\geq s_1$, since $S''$ is a final odd segment of $S_1$ and $S_1$ satisfies $\ASC$.
It follows that the weight sequence $(s',s',s_2,\ldots,s_m)$ is again weakly $U$--shaped, and the theorem implies $\tes(a'_2,a_3,\ldots,a_n)\sim\langle s', s',s_2,\ldots,s_m\rangle \sim \langle s_2,\ldots,s_m\rangle$.

(ii)\ The proof is the same as the proof of (i) with minor changes, such as replacing ``weakly $U$--shaped'' with ``non--increasing''.
\enddemo

If $g=g_1\oplus\ldots \oplus g_n$ where each $g_i$ is a stack or two--ended stack, then Lemma 2.10 and Theorem 2.11 give an equivalence of $g$ with a game $g'$ of the form $\langle s_1,\ldots,s_m\rangle \oplus \ev(x)$.
It is easy to compute the value of $g'$.
If $m>0$, then the largest of the $s_i$'s is an optimal initial move in $g'$, and Corollary 2.16 implies that the corresponding move in $g$ is optimal.
If $m=0$, note that $g$ is actually $\bigoplus_1^n\st (\ux_i)$, where each $\ux_i$ is an $\ev$--sequence.
In this case, it is easily seen that any initial move in $g$ is optimal.  If Player One chooses an initial move as just indicated, it is optimal for both players to continue moving in this slice or $ev-$sequence until it is exhausted, and thus we have an optimal scenario for a complete play of the game.

So we have solved games of this type modulo the computation of the partitions as in Lemma 2.10.
In the course of play of the game $g$ it may be necessary for a player to perform this computation before each move.
For example, if $m=0$ any initial move is optimal, but there may be optimal responses other than continuing in this stack. 
Thus situations where recalculation is necessary could arise under optimal play, even in this seemingly simple case.

We now present an algorithm for computing the desired partition in linear time in 2.17--2.20 below.

\medskip\noindent
{\bf 2.17\ Description of the algorithm}

I.\ Two--ended stacks.
The algorithm consists of three parts.
In part A we split off zero or more slices from the left of the given sequence, in part B we split off zero or more slices from the right, and in part C we deal with what remains in the middle.

\noindent
A.1.\ Given the sequence $(a_1,\ldots,a_n)$, seek the smallest even number $p\leq n$ such that $\sum_1^p (-1)^{i-1} a_i>0$.
If no such $p$ exists, go on to part B.
If $p$ exists, go to step 2.

2.\ Find the largest odd number $j<p$ such that $\sum_1^j (-1)^{i-1} a_i$ is minimal.
Then $(a_1,\ldots,a_j)$ is a slice (and is the longest initial segment that is a slice).
Save this slice, replace $(a_1,\ldots,a_n)$ with $(a_{j+1},\ldots,a_n)$, and go back to step 1.

\noindent
B.\ This is the same as part A with left and right reversed.
So when we are done, we will have split off slices from the right, and the reverse of the remaining sequence will satisfy $\ASC$.
Then go on to part C.

\noindent
C.\ Now our input is a sequence $\ub=(b_1,\ldots,b_q)$ such that both $\ub$ and its reverse satisfy $\ASC$.
If $q$ is odd, $\ub$ is the last slice and we are done.
If $q$ is even, $\ub=S' S''$, where $S'$ and $S''$ are slices of the same weight.
To find this decomposition, find the largest odd number $j<q$ such that $\sum_1^j (-1)^{i-1} b_j$ is minimal, and let $S'=(b_1,\ldots,b_j)$, similarly to step A.2 and part of the proof of corollary 2.16.

\noindent
II.\ Stacks.
Here the algorithm consists of two parts.
Part A splits off zero or more slices from the left and part B deals with what remains.

A.\ This is identical to I.A.

B.\ Our input is a sequence $\ub=(b_1,\ldots,b_q)$ which satisfies $\ASC$.
If $q$ is even, let $X=\ub$, which is an $\ev$--sequence.
If $q$ is odd, $\ub$ is $S$ or $SX$ where $S$ is a slice and $X$ is an $\ev$--sequence.
To find this, find the largest odd number $j\leq q$ such that $\sum_1^j (-1)^{i-1} b_i$ is minimal and let $S=(b_1,\ldots,b_j)$.

\medskip\noindent
{\bf 2.18\ Justification of the algorithm}.

I.\ We know from Lemma 2.10(i) that $(a_1,\ldots,a_n)=S_1\ldots S_m$, where the $S_i$'s are slices of weight $s_i$ and $(s_1,\ldots,s_m)$ is $U$--shaped.
We will use this as a model to analyze what happens in the algorithm and show that the algorithm produces slices with the same weights.

In step A.1, if the term $a_p$ lies in $S_j$ with $j$ even, then $\sum_1^p (-1)^{i-1} a_p\leq \sum_1^j (-1)^{i-1} s_i$.
And if $a_p$ lies in $S_j$ with $j$ odd, then $j>1$ and $\sum_1^p (-1)^{i-1} a_p\leq\sum_1^{j-1} (-1)^{i-1} s_i$.
It follows that if $s_1>s_2$, then $a_p$ lies in $S_2$.
And if $s_1\leq s_2$ (or if $m=1$), then $p$ does not exist.

If $p$ exists, then the fact that the reverse of $S_1$ satisfies $\ASC$ implies that $s_1\leq\sum_1^j (-1)^{i-1} a_i$ for odd numbers $j$ less than the size of $S_1$.
And the fact that $S_2$ satisfies $\ASC$ implies that $s_1\leq \sum_1^j (-1)^{i-1} a_i$ for odd numbers $j$ greater than the length of $S_1$ (and $j<p$).
So $s_1$ is the minimum obtained in step A.2 and $j$ is at least the size of $S_1$.
Thus $(a_1,\ldots,a_j)$ is $S_1\ub$, where $\ub$ is an even initial segment of $S_2$ which is an $\ev$--sequence of weight 0.
This implies that the reverse of $\ub$ is also an $\ev$--sequence, and hence $(a_1,\ldots,a_j)$ is indeed a slice of weight $s_1$.
Also, if $S'_2$ is the result of removing $\ub$ from $S_2$, then $S'_2$ is a slice of weight $s_2$.
Going forward from here we replace $S_2$ by $S'_2$ and change notation.
In other words, in the new notation, $S_2\ldots S_m$ is a model for $(a_{j+1},\ldots,a_n)$.

Now it is clear that part A of the algorithm will split off slices of weights (in order) $s_1,\ldots,s_j$ where $s_j>s_{j+1}\leq s_{j+2}$ (or $j+1=m$).
Moreover, by modifying the slices $S_1,\ldots,S_{j+1}$ in our model decomposition as indicated above, we may assume that the input into part B is $S_{j+1}\ldots S_m$.

Exactly the same analysis applies to part B.
So we see that the sequence inputted into part C will consist of exactly the one or two (possibly modified) $S_i$'s of minimal weight.
If only one of the $S_i$'s has minimal weight, the $q$ in part C will be odd.
If two $S_i$'s have minimal weight, we need only observe the following:\ \ If $\ub=(b_1,\ldots,b_q)$, $q$ is even, both $\ub$ and its reverse satisfy $\ASC$, and $\ub=S'S''$ where $S'$ and $S''$ are slices, then the weight of $S'$ (and $S''$) must be the minimum of $\sum_1^j(-1)^{i-1} b_j$ for odd $j<q$, and for any odd $j$ which achieves the minimum, both $(b_1,\ldots,b_j)$ and $(b_{j+1},\ldots,b_q)$ are slices.

\noindent
II.\ We now start with $(a_1,\ldots,a_n)=S_1\ldots S_m$ or $S_1\ldots S_m X$, as in Lemma 2.10 (ii), and use this to analyze what happens in the algorithm.
If $m=0$, then obviously the algorithm produces $X$ (and no slices).
If $m\geq 2$, then the argument in $I$ above shows that the first $m-1$ run--throughs of A.1 and A.2 will produce slices of weights $s_1,\ldots,s_{m-1}$ in order and leave a remaining sequence consisting of $S_m$ or $S_m X$ (with a possibly modified $S_m$).
This is where we start if $m=1$, so we now assume $(a_1,\ldots,a_n)=S_m$ or $S_m X$.
There are two possibilities:\ \ Either step A.1 finds a number $p$ or it doesn't.

If $p$ exists, then the argument in $I$ above shows that we split off a slice of weight $s_m$, and what remains is $X$ or the result of removing from $X$ an initial $\ev$--sequence of weight 0.
So we are done in this case.
If $X$ has weight 0, it is possible that nothing will be left after removing the last slice.

If $p$ does not exist (a case which includes the possibility that $X$ is missing), then the input into part B will be odd.
Then arguments, like those in I show that we will split off a slice $S'$ with weight $s_m$.
Moreover $S'$ will either be $S_m$ or $S_m$ followed by an $\ev$--sequence of weight 0.
And what remains will either be $X$ or $X$ with an initial even segment, an $\ev$--sequence of weight 0, removed.
If $X$ has weight 0, all of $X$ will be part of the last slice.

\medskip\noindent
{\bf 2.19.\ How long does the algorithm take?}

Note that if step A.1 produces a number $p$, then the number of steps required to find $p$ is only $O(p)$.
(To calculate $\sum_1^{q+2} (-1)^{i-1} a_i$ we just add two more terms to the previously calculated --- and temporarily saved --- sum, $\sum_1^q (-1)^{i-1} a_i$.)
Then the number of additional steps needed to complete step A.2 is also only $O(p)$.
Now suppose $k$ slices are split off by part A and the sizes of the first $k+1$ slices of our model are $\ell_1,\ldots,\ell_{k+1}$.
Then the first slice requires $O(\ell_1+\ell_2)$ steps, the second $O(\ell_2+\ell_3)$ steps, etc.
After the first $k$ slices have been removed, the last run--through of step A.1, which confirms that part A is finished, requires $O(n-\ell_1-\ldots-\ell_k)$ steps.
Thus the total number of steps for part A is bounded by a constant times $(\ell_1+\ell_2)+\ldots+(\ell_k+\ell_{k+1})+n-\ell_1-\ldots-\ell_k\leq 2n$.
Similarly, part I.B requires only $O(n)$ steps, and obviously part C also requires only $O(n)$ steps.
(Note that only $q-1$ comparisons are needed to find the minimum of $q$ numbers.)
Similar arguments show that version II of the algorithm also requires only $O(n)$ steps.

We are not providing an estimate for a number $K$ such that the number of steps required is at most $Kn$, but it is clear that $K$ is not unduly large.

\medskip\noindent
{\bf 2.20.\ Is the algorithm computationally effective and stable?}

There are two issues here.
One is that some of the comparisons or sign determinations may be very close calls, so that a small error could cause problems.
The other is that if the algorithm has any applications outside of game theory, the data inputted may be only approximate.
So we would like to know that small errors in the input produce only small errors in the output.
Without claiming to be definitive, we provide three comments pointing towards a positive answer to the question.

1.\ In the case of games that people might actually play, it is reasonable to suppose that the weights will be rational.
Then we can find a common denominator and do exact integer arithmetic with the numerators.
If some of the numerators are large, it may be necessary to use multiple precision.
So the estimate for running time should take this into account.

2.\ Another approach would be to round all the data to a certain number of decimal places and use enough multiple precision so that no additional rounding is necessary.
Aside from the error caused by the initial round-off, this amounts to the same thing as above with the common denominator being a power of 10.

3.\ We already have a metric which can be used to compare the outputs of the algorithm for two (close) inputs, $\ua=(a_1,\ldots,a_n)$ and $\ua'=(a'_1,\ldots,a'_n)$, namely the distance function on games introduced in \S1.
Suppose the outputs of version I of the algorithm are $(s_1,\ldots,s_m)$ for the input $\ua$ and $(s'_1,\ldots,s'_{m'})$ for $\ua'$.
And suppose the outputs of version II are $(s_1,\ldots,s_k;x)$ and $(s'_1,\ldots,s'_{k'},x')$.
Then it follows from Lemma 1.3 and Theorem 2.11 that the outputs will be close in our metric if $\ua$ and $\ua'$ are close in the $\ell^1$--sense.
More precisely,
$$
\sI (\langle s_1,\ldots,s_m\rangle, \langle s'_1,\ldots,s'_{m'}\rangle)\leq \|\ua-\ua'\|_1.\tag2.1
$$
$$
\sI(\langle s_1,\ldots,s_k\rangle, \langle s'_1,\ldots,s'_{k'}\rangle) + |x-x'|\leq \|\ua-\ua'\|_1.\tag2.2
$$
(It will be shown in \S3 below that 
$$
\gathered
\sI(\langle s_1,\ldots,s_k\rangle\oplus\ev(x), \langle s'_1,\ldots,s'_{k'}\rangle\oplus\ev(x'))=\\
\sI(\langle s_1,\ldots,s_k\rangle,\ \langle s'_1,\ldots,s'_{k'}\rangle)+|x-x'|.)
\endgathered
$$
For game theory this seems to be an adequate stability result.
Closeness in our metric does not rule out the possibility that $m$ and $m'$ could be very different, but this seems appropriate.
For hypothetical applications of the algorithm outside of game theory, it seems to be at least a reasonable hope that (2.1) and (2.2) will be adequate.
(The metric for games of the type $\langle s_1,\ldots,s_m\rangle$ is described more explicitly in \S3 below).

{\bf \S3.\ Some general theory}.

All of the formal results of this section are stated for graph games, but many of the arguments work also for general games with predetermined parity.
But if we worked in the more general class, then the concepts of distance and equivalence would change:\ \ In the relation $\sI(g_1,g_2)=\sup|\val(g_1\oplus h)-\val (g_2\oplus h)|$, $h$ would range over a wider class of games and $\oplus$ would be replaced by $\vee$.
However, as indicated in Remark 3.26 below, for an even graph game $g$, $\sI(g,\langle \ \rangle$) is the same regardless of whether $h$ is required to be a graph game or allowed to be a general game (even without parity).
In fact to calculate $\sI(g,\langle \ \rangle$) we ``merely'' have to calculate $\val(g\oplus \langle y\rangle$) for $y < -|g|$.
It follows that any invertible graph game is still invertible within the class of general games, and obviously $\sI(g,h)=\sI(g\oplus h,\langle \ \rangle)$ when $g$ is invertible.
The reader who is so inclined can keep these considerations in mind when reading this section.

\proclaim{Proposition 3.1}(i)\ If $g$ and $h$ are even graph games, then $\val(g\oplus h)\leq \val(g)+\val(h)$.

(ii)\ If $g$ is an odd graph game and $h$ an even graph game, then $\val(g\oplus h)\geq\val(g)-\val(h)$.
\endproclaim

\demo{Proof}(i)\ We use induction on the size of $g\oplus h$, the case of size 0 being trivial.
If Player One makes an initial move $v$ in $g$, then Player Two can respond with a move $w$ in $g\backslash v$ which is optimal for the game $g\backslash v$.
Then Player One is looking at $g\backslash \{v,w\}\oplus h$ , which by induction has value at most $\val(g\backslash \{v,w\})+\val(h)$.
Thus the outcome is at most $\val(g)+\val(h)$.
A similar argument applies if Player One chooses an initial move in $h$.

(ii)\ Player One can begin with a move $v$ which is optimal for $g$.
Then apply part (i) to $(g\backslash v)\oplus h$.
So the outcome is $\wt(v)-\val((g\backslash v)\oplus h)\geq \wt(v)-\val(g\backslash h)-\val(h)=\val(g)-\val(h)$.
\enddemo

\proclaim{Proposition 3.2}For any graph game $g$, $\val(g\oplus g)\leq 0$.
\endproclaim

\demo{Proof}Player Two can use the well known symmetric strategy.
\enddemo

A graph game $g$ will be called {\it invertible} if there is a graph game $h$ such that $[g]+[h]=[g\oplus h]=\bold 0$.
Let $\sG=\{ [g]\colon g$ is an invertible graph game$\}$.
So $\sG$ is the largest group containing $\bold 0$ and contained in the semigroup $\sS$.

\proclaim{Theorem 3.3}If $g$ is an invertible graph game, then $[g]+[g]=\bold 0$.
\endproclaim

\demo{Proof}Suppose $h$ is a graph game such that $g\oplus h\sim \langle\ \rangle$.
Then $g$ and $h$ have the same parity by Proposition 1.4.
Let $k$ be any graph game.
Then
$$
g\oplus k\sim g\oplus (g\oplus h)\oplus k\sim (g\oplus g)\oplus (h\oplus k).\tag3.1
$$
If $h\oplus k$ is even, then Propositions 3.1 and 3.2 imply
$$
\val (g\oplus k)\leq \val(g\oplus g)+\val (h\oplus k)\leq \val (h\oplus k).
$$
Since $g\oplus k$ is also even and the hypothesis on $(g,h)$ is symmetric, we also have $\val (h\oplus k)\leq \val (g\oplus k)$; i.e., $\val(g\oplus k)=\val (h\oplus k)$.
If $h\oplus k$ is odd, then (3.1) implies 
$$
\val (g\oplus k)\geq \val(h\oplus k)-\val (g\oplus g)\geq \val(h\oplus k).
$$
Again symmetry implies $\val(g\oplus k)=\val (h\oplus k)$.
So we have proved $g\sim h$.
\enddemo

\proclaim{Theorem 3.4}If $g$ is a graph game, then a necessary condition for $g$ to be invertible is that $\val(g\oplus g)=0$, and a sufficient condition is that $\val(g'\oplus g')=0$ for all games $g'$ that can arise from legal play in $g$.
\endproclaim

\demo{Proof}The necessity follows immediately from Theorem 3.3.
For the sufficiency we need to prove that $\val(g\oplus g\oplus h)=\val(h)$ for all graph games $h$, and we use induction on the size of $g\oplus h$, the case of size 0 being obvious.
The induction hypothesis implies that for any vertex $v$ of $h$, the outcomes of the two $v$--moves are the same:
$$
\wt(v)-\val (g\oplus g\oplus (h\backslash v))=\wt(v)-\val (h\backslash v).
$$
If $h$ is not the empty game, this implies $\val (g\oplus g\oplus h)\geq \val(h)$; and if $h$ is empty, the same inequality follows from $\val(g\oplus g)=0$.
To show that $\val (g\oplus g\oplus h)\leq \val(h)$, we need only show that if $w$ is a vertex of $g$ which is a legal initial move, then the $w$--move in one of the copies of $g$ doesn't lead to a greater outcome than $\val(h)$.
But Player Two can respond with the $w$--move in the other copy of $g$.
Then Player One is looking at $(g\backslash w)\oplus (g\backslash w)\oplus h$, and the induction hypothesis together with the fact that $g\backslash w$ satisfies the same hypothesis as $g$ implies the result.
\enddemo

\proclaim{Theorem 3.4$'$}A sufficient condition for the graph game $g$ to be invertible is that $\val(g\oplus g)=0$ and every legal initial move $w$ in $g$ is metrically dominated by a legal initial move $v$ such that $g\backslash v$ is invertible.
\endproclaim

\demo{Proof}The proof is very similar, but now, to prove $\val(g\oplus g\oplus h)=\val (h)$, we use induction just on the size of $h$.
The only difference in the proof is the step where we show that no legal initial move in one of the copies of $g$ can lead to a greater outcome in $g\oplus g\oplus h$ than $\val(h)$.
But by hypothesis we may assume this move is a vertex $v$ such that $g\backslash v$ is invertible.
Then Player Two can respond with the $v$--move in the other copy of $g$.
Since $\val (g\backslash v\oplus g\backslash v\oplus h)=\val(h)$, the result follows.
\enddemo

\example{Remark}Theorem 3.4$'$ is stronger than Theorem 3.4 in the sense that any game which satisfies the sufficient condition of Theorem 3.4 also satisfies the hypothesis of Theorem 3.4$'$.
However, we need Theorem 3.4 to prove this assertion.
\endexample

\example{Example 3.5}The game $g_1=\path(1,2,3)$ is invertible by Theorem 3.4$'$ but does not satisfy the sufficient condition of Theorem 3.4.
We do not know of any invertible graph games which fail to satisfy the sufficient condition of Theorem 3.4$'$.
The game $g_2=\path (0,1,0,1)$ satisfies $\val (g_2\oplus g_2)=0$ but is not invertible.
To see that it is not invertible, compute that $\val(g_2\oplus g_2\oplus \langle y\rangle)=y+2$ if $y<<0$.
(Theorem 4.3 below can simplify the verification that $g_1$ satisfies the hypothesis of Theorem 3.4$'$ as well as the verification of the assertions about $g_2$.)
\endexample

\proclaim{Proposition 3.6}If $g$ is any graph game and $h$ is an invertible even graph game then
$$
\val(g)-\val(h)\leq \val(g\oplus h)\leq \val(g)+\val (h).
$$
\endproclaim

\demo{Proof}If $g$ is even, the righthand  inequality follows from Proposition 3.1.
Also, since $g\sim g\oplus (h\oplus h)\sim (g\oplus h)\oplus h$, Proposition 3.1 yields $\val(g)\leq \val(g\oplus h)+\val (h)$, which implies the lefthand inequality.
A similar argument works if $g$ is odd.
\enddemo

\proclaim{Corollary 3.7}If $h$ is an invertible even graph game, then $\val(h)\geq 0$.
\endproclaim

\proclaim{Corollary 3.8}If $g_1$ and $g_2$ are invertible graph games of the same parity, then $\sI (g_1,g_2)=\val (g_1\oplus g_2)$.
\endproclaim

\demo{Proof}Since $\val(g_1\oplus g_2)-\val (g_2\oplus g_2)=\val (g_1\oplus g_2)$, then $\sI(g_1,g_2)\geq \val(g_1\oplus g_2)$.
If $k$ is any graph game, then $|\val (g_1\oplus k)-\val (g_2\oplus k)|=$$|\val((g_1\oplus g_2)\oplus (g_2\oplus k))-\val (g_2\oplus k)|\leq \val(g_1\oplus g_2)$, by the Proposition.
\enddemo

\proclaim{Proposition 3.9}If $g$ is an invertible odd graph game and $h$ is any odd graph game, then $\val (g\oplus h)\geq \val (g)-\val (h)$.
\endproclaim

\demo{Proof}Since $h\sim g \oplus (g\oplus h)$, Proposition 3.1 implies $\val(h)\geq \val(g)-\val (g\oplus h)$, which yields the result.
\enddemo

\example{Remark 3.10}Consider the inequalities
$$
\align
&\val (g\oplus h)\leq \val(g) +\val (h),\text{ and}\tag3.2\\
&\val(g\oplus h)\geq \val(g)-\val (h).\tag3.3
\endalign
$$
Propositions 3.1, 3.6, and 3.9 imply these inequalities under various hypotheses.
The hypotheses involve only the parities of $g$ and $h$ or invertibility of $g$ or $h$.
These propositions include all cases where (3.2) or (3.3) follows from a hypothesis of this type.
Thus if $g$ and $h$ are odd, (3.2) can fail even if both are invertible (take $g=h=\langle -1\rangle$).
And if $g$ and $h$ have opposite parity (3.2) can fail if the odd one is invertible (take $g=\langle 0\rangle$ and $h$ in $\ev(1)$).
Also, if $g$ is even and $h$ odd, (3.3) can fail even if both are invertible (take $g=\langle \ \rangle$ and $h=\langle -1\rangle$).
If both are even, (3.3) can fail if $g$ is invertible (take $g=\langle \ \rangle$ and $h$ in $\ev(1)$).
And if both are odd, (3.3) can fail if $h$ is invertible (take $g$ of the form $\langle 0\rangle \oplus \ev(1)$ and $h=\langle 0\rangle$).
\endexample

\example{3.11.\ The completions of $\sG_0$ and $\sS_0$}
Let $\mu$ be Lebesgue measure on the real line, $\bR$, and let $\sM$ be its measure algebra.
So the elements of $\sM$ are equivalence classes of measurable subsets of $\bR$, with two sets equivalent if they differ by a set of measure zero; and (finite or countable) operations on measurable subsets of $\bR$ carry over to $\sM$.
For $A$ and $B$ in $\sM$, let $A+B$ denote the symmetric difference, $(A\backslash B)\cup (B\backslash A)$.
Then $\sM$ becomes an abelian group and is a vector space over $\bF_2$, the field with two elements.
Note that $\sG_0$ is also a vector space over $\bF_2$ and $\{ [x]\colon x\in\bR\}$ is a basis.
Thus we can define a homomorphism $\theta\colon \sG_0\to\sM$ by $\theta ([x])=(-\infty,x)$.
\endexample

If $\bg=[x_1,\ldots,x_n]$ is an even element of $\sG_0$, with $x_1 < x_2<\ldots < x_n$, then $\theta(\bg)=(x_1,x_2)\cup\ldots\cup (x_{n-1},x_n)$.
Therefore $\mu(\theta(\bg))=\val(\bg)$.
Now the distance function on $\sM$ is given by $\dist(A,B)=\mu(A+B)$.
So it follows from Corollary 3.8 that for $\bg_{\bold 1}$ and $\bg_{\bold 2}$ of the same parity in $\sG_0$, $\dist(\theta (\bg_{\bold 1}), \theta(\bg_{\bold 2}))=\sI(\bg_{\bold 1},\bg_{\bold 2})$.
It is clear then that $\theta$ extends to an isometry, still denoted by $\theta$, from the completion $\overline{\sG}_0$ of $\sG_0$ onto the closure of $\theta(\sG_0)$.
Now elementary measure theory shows that the closure of $\{\theta (\bg)\colon \bg\in\sG_0$, $\bg$ even$\}$ is $\{A\in\sM\colon \mu(A) < \infty\}$.
Every odd element of $\sG_0$ is of the form $[0]+\bg$, $\bg$ even.
Therefore the closure of $\{\theta (\bg)\colon \bg\in \sM,\bg\text{ odd}\}$ is $\{A\in\sM\colon\mu (A+(-\infty,0)) < \infty\}$.

In order to describe $\val(\bg)$ in terms of Lebesgue measure when $\bg$ is odd, we define a function $v$ on $\{A\in \sM\colon \mu (A+(-\infty,0))<\infty\}$:
$$
\text{For }c\in\bR, v(A)=c+\mu(A\backslash (-\infty,c))-\mu((-\infty,c)\backslash A).\tag3.4
$$
It is easy to check that the result of (3.4) is independent of the choice of $c$.
(If $c_1<c_2$, just look separately at $A\cap (c_2,\infty), A\cap (c_1,c_2)$, and $A\cap (-\infty,c_1)$.)
Now if $\bg=[x_1,\ldots,x_n]$ is an odd element of $\sG_0$, with $x_1<x_2\ldots <x_n$, we can confirm that $\val (\bg)=v(\theta (\bg))$ by choosing $c\leq x_1$.

\proclaim{Proposition 3.12}If $g_1$ and $g_2$ are invertible graph games of the same parity and $x_1,x_2\geq 0$, then
$$
\sI (g_1\oplus\ev(x_1), g_2\oplus\ev(x_2))=\sI(g_1,g_2)+|x_1-x_2|.
$$
\endproclaim

\demo{Proof}It is obvious that 
$\sI(g_1\oplus\ev(x_1), g_2\oplus\ev (x_2))\leq\sI(g_1,g_2) + |x_1-x_2|$.
To show the reverse, we may assume without loss of generality that $x_1\leq x_2$.
Then
$$
\gathered
\sI(g_1\oplus \ev(x_1),g_2\oplus\ev(x_2))\geq |\val (g_1\oplus\ev(x_1)\oplus g_2)-\val(g_2\oplus\ev (x_2)\oplus g_2)|\\
=|\val (g_1\oplus g_2)-x_1+x_2|=\sI (g_1,g_2)+|x_1-x_2|.
\endgathered
$$
\enddemo

\example{Remark}The game called $g_2$ in Example 3.5 is not of the form $g\oplus\ev(x)$, $g$ invertible.
\endexample

Now it is clear that $\overline{\sS_0}$, the completion of $\sS_0$, can be identified with $\overline{\sG_0}\oplus [0,\infty)$.

Motivated by the fact that $\sG_0$ and $\overline{\sG_0}$ have been identified with Boolean rings (though we have no {\it a priori} concept of multiplication on $\sG_0$), we will define two relations on $\sG$, the group of all equivalence classes of invertible graph games.
The usefulness of these relations is not so clear because neither $\sG$ nor its completion is even a lattice.

\definition{Definition 3.13}(i)\ If $\bg_{\bold 1},\bg_{\bold 2}\in\sG$, we say $\bg_{\bold 1}\geq \bg_{\bold 2}$ if either $\bg_{\bold 1}$ is odd or $\bg_{\bold 2}$ is even and
$$
\val(\bg_1+\bg_{\bold 2})=\val (\bg_1)-\val(\bg_{\bold 2})
$$

(ii)\ If $\bg_{\bold 1},\ldots,\bg_{\bn}\in\sG$, we say that $\{\bg_{\bold 1},\ldots,\bg_{\bn}\}$ is {\it independent} if at most one of the $\bg_{\bi}$'s is odd and 
$$
\val(\bg_{\bold 1}+\ldots+\bg_{\bn})=\val (\bg_{\bold 1})+\ldots+\val (\bg_{\bn}).
$$
If $n=2$, we also say that $\bg_{\bold 1}$ and $\bg_{\bold 2}$ are independent.

Note that if we replace = by $\geq$ in (i) or = by $\leq$ in (ii), we get a relation that is always true.
Also $\bg_{\bold 1}\geq \bg_{\bold 2}$ if and only if $\bg_{\bold 1}+\bg_{\bold 2}$ and $\bg_{\bold 2}$ are independent.
And it is easy to see that if the $\bg_{\bold i}$'s are in $\sG_0$, then $\bg_{\bold 1}\geq \bg_{\bold 2}$ if and only if $\theta (\bg_{\bold 1})\supset \theta
(\bg_{\bold 2})$ and $\{\bg_{\bold 1},\ldots,\bg_{\bn}\}$ is independent if and only if $\theta(\bg_{\bold 1}),\ldots,\theta(\bg_{\bn})$ are mutually disjoint.
\enddefinition

\proclaim{Proposition 3.14}Assume $\bg_{\bold 1},\ldots,\bg_{\bn}\in\sG$ and the set $\{1,\ldots,n\}$ is partitioned into disjoint subsets $J_1,\ldots,J_k$.
Let $\bh_{\bj}=\sum_{i\in J_j} \bg_{\bi}$.
Then $\{\bg_{\bold 1},\ldots,\bg_{\bn}\}$ is independent if and only if each set $\{\bg_{\bi}\colon i\in J_j\}$ is independent and $\{\bh_{\bold 1},\ldots,\bh_{\bk}\}$ is independent. 
\endproclaim

\demo{Proof}Since $\bg_{\bold 1}+\ldots+\bg_{\bn}=\bh_{\bold 1}+\ldots+\bh_{\bk}$, we have 
$$
\gathered
\val (\bg_{\bold 1}+\ldots+\bg_{\bn})=\val (\bh_{\bold 1}+\ldots+\bh_{\bk})\leq\\
\val(\bh_{\bold 1})+\ldots+\val (\bh_{\bk})\leq\sum_1^k \sum_{i\in J_j}\val (\bg_{\bi})=\\
\val (\bg_{\bold 1})+\ldots+\val (\bg_{\bn}).
\endgathered
$$
Note that under either of the conditions being proved equivalent, all of the parity requirements are satisfied.
If $\{\bg_{\bold 1},\ldots,\bg_{\bn}\}$ is independent, then the outside parts of this relation string are equal, and therefore all of the inequalities are equalities, which implies the second condition.
And if the second condition is satisfied, all of the inequalities are equalities.
\enddemo

\proclaim{Corollary 3.15}Any subset of an independent set is independent, and the sums of two disjoint subsets
of an independent set are independent of one another.
\endproclaim

\proclaim{Corollary 3.16}If $\{\bg_{\bold 1},\ldots,\bg_{\bn}\}$ is independent and $\bg_{\bi}\geq \bh_{\bi},\forall i$, then $\{\bh_{\bold 1},\ldots,\bh_{\bn}\}$ is independent.
\endproclaim

\demo{Proof}Note that $\bg_{\bi}=\bh_{\bi}+ (\bg_{\bi}+\bh_{\bi})$ and $\bh_{\bi}$ is independent of $\bg_{\bi}+\bh_{\bi}$.
\enddemo

\proclaim{Proposition 3.17}The relation $\geq$ is a proper order relation on $\sG$.
\endproclaim

\demo{Proof}It is obvious that $\bg\geq\bg$.
If $\bg\geq\bh$ and $\bh\geq\bg$, then $\val(\bg+\bh)=\val (\bg)-\val(\bh)=\val(\bh)-\val(\bg)$, and also $\bg$ and $\bh$ have the same parity.
So $\val(\bg+\bh)=0$ and $\bg=\bh$ by Corollary 3.8.
If $\bg\geq\bh$ and $\bh\geq\bk$, then $\bg+\bh$ is independent of $\bh$ and $\bh+\bk$ is independent of $\bk$.
Since $\bh=\bk+(\bh+\bk)$, this implies $\{\bg+\bh,\bh+\bk,\bk\}$ is independent by Proposition 3.14.
It follows that $\bk$ is independent of $(\bg+\bh)+(\bh+\bk)=\bg+\bk$, whence $\bg\geq\bk$.

The next example, combined with Proposition 3.14, explains why we chose the term ``independent'' rather than the intuitively stronger term ``mutually disjoint'' in Definition 3.13.
\enddemo

\example{Example 3.18}Let $g=\cyc (0,1,0,1,1,0), h_1=\langle 0,1\rangle$, and $h_2=\langle -1,0\rangle$.
Then $[g]$ is independent of both $[h_1]$ and $[h_2]$ but not independent of $[h_1\oplus h_2]$.
Also $[h_1]$ and $[h_2]$ are independent.
(It will be shown in the next section that all unbroached cycle games are invertible.)
We leave it to the reader to verify that $\val(g)=\val(h_1)=\val(h_2)=1$, $\val(g\oplus h_1)=\val(g\oplus h_2)=\val(h_1\oplus h_2)=2$, 
and $\val(g\oplus (h_1\oplus h_2))=1$.
Actually $[g]\leq [h_1\oplus h_2]$.
\endexample

\example{Example 3.19}$\sG$ is not a lattice.
We exhibit even invertible games $g_1,g_2,h_1,h_2$ such that $[g_1],[g_2]\geq [h_1],[h_2]$ but there does not exist $k$ with $[g_1],[g_2]\geq [k]\geq [h_1],[h_2]$.
Thus $\sG$ does not satisfy the Riesz interpolation property, and {\it a fortiori} it is neither a meet semilattice nor a join semilattice.
Again let $h_1=\langle 0,1\rangle$ and $h_2=\langle -1,0\rangle$, but now take $g_1=\langle -1,1\rangle$ and $g_2=$ $\cyc(0,1,0,1,1,0,1,0)$.
Of course, $[g_1]$ is the least upper bound of $[h_1]$ and $[h_2]$ in $\sG_0$, but also $[g_2]\geq [h_1], [h_2]$, and $\val (g_1)=\val(g_2)=2$.
If either $[g_1]\geq [g_2]$ or $[g_2]\geq [g_1]$, then we would have $g_1\sim g_2$, but in fact $\val(g_1\oplus g_2)=2$, so this is not the case.
If $[k]$ exists as above then $[k]\leq [g_1], [g_2]$ and $[k]\neq [g_1], [k]\neq [g_2]$.
This implies that $k$ is even and $\val(k)<2$.
Then $[k]\geq [h_i]$ implies $\val(k\oplus h_i)=\val(k)-\val (h_i) < 1$, whence $\val(h_1\oplus h_2)=\val ((h_1\oplus k)\oplus (k\oplus h_2))\leq \val(h_1\oplus k)+\val (k\oplus h_2) < 2$, a contradiction.
The concepts and arguments carry over to $\overline{\sG}$, the completion of $\sG$, and hence $k$ does not exist in $\overline{\sG}$ either.
\endexample

The next example exhibits some elements of $\overline{\sG}$ which are variants of the familiar I--cut--and--you--choose game.

\example{Example 3.20}For $x>0$, we define two games $\gamma(x)$ and $\sigma(x)$, which are not graph games, but the arguments of this section still work for them.
In both games the initial move is a number $a$ in $(0,x)$ and has weight 0.
In $\gamma(x)$, Player Two is then faced with the graph game $\langle a,x-a\rangle$, and in $\sigma(x)$ Player Two is faced with $\langle -a,-(x-a)\rangle$.
So $\sigma(x)$ can be regarded as the I--cut--and--you--serve game.
Both $\gamma(x)$ and $\sigma(x)$ are invertible odd games of value 0, but they are not as similar as they may seem.
In fact, $\gamma(x)\sim\langle 0\rangle$, $\forall x>0$, and $\sI(\sigma(x),\langle 0\rangle)=x$.
So far as we know, $\sigma(x)$ is not equivalent to any graph game, but $\sigma(x)$ is the limit of a sequence of graph games.
Let $g_n(x)=\path (0,x/n,0,x/n,\ldots,0)$, where there are $n$ weights of $x/n$ and $n+1$ weights of 0.
It follows from the results of \S2, that if $v$ is a vertex of weight $x/n$ in $g_n(x)$, then $g_n(x)\backslash v\sim \langle -(k/n)x, -(\ell/n)x\rangle$, where $k$ and $\ell$ are non--negative integers such that $k+\ell=n-1$.
It follows from Theorem 4.3 below that each 0--move in $g_n(x)$ is metrically dominated by an $x/n$--move, or this can be computed directly.
These facts show that $g_n(x)$ is invertible and $g_n(x)\to\sigma(x)$.
\endexample

\medskip\noindent
{\bf 3.21\ Action of the $(ax+b)$--group}.

Consider the transformations $\varphi_{a,b}\colon\bR\to\bR$ defined by $\varphi_{a,b}(x)=ax+b$, where $a>0,\ b\in\bR$.
For any graph game $g$, let $\varphi_{a,b}(g)$ be obtained from $g$ by replacing each weight $w$ with $\varphi_{a,b}(w)$.
If $a=1$, then $\varphi_{a,b}(g)=\tau_b(g)$, as previously defined.
It is clear that the strategy of $\varphi_{a,b}(g)$ is identical to that of $g$ and $\val(\varphi_{a,b}(g))=\varphi_{a,b}(\val(g))$ if $g$ is odd, and $\val(\varphi_{a,b}(g))=$  $a\,\val(g)$ if $g$ is even.
We see then that $\sI(\varphi_{a,b}(g_1)$, $\varphi_{a,b}(g_2))=$  $a\,d(g_1,g_2)$, and in particular $g_1\sim g_2$ implies $\varphi_{a,b}(g_1)\sim \varphi_{a,b}(g_2)$.
(In comparing $\val(\varphi_{a,b} (g_1)\oplus h)$ with $\val(\varphi_{a,b}(g_2)\oplus h)$, note that $\varphi_{a,b}(g_i)\oplus h=\varphi_{a,b}(g_i\oplus\varphi_{a,b}^{-1}(h))$.)
Therefore we may write $\varphi_{a,b} ([g])=[\varphi_{a,b}(g)]$, and we have a jointly continuous action of the $(ax+b)$--group on $\sS$, which extends to an action on the completion $\overline{\sS}$.
For $[g]$ in $\sG_0$, $\theta (\varphi_{a,b} (g))=\Phi_{a,b}(\theta(g))$, where $\Phi_{a,b}(A)=\{\varphi_{a,b}(x)\colon x\in A\}$, $A\subset\bR$; and of course $\Phi_{a,b}$ operates also on $\sM$.

It is not hard to see that if $g$ is an even graph game and $h$ is any graph game, then $\val(\tau_\lambda(g)\oplus h)\geq \val(g)+\val(h)$ for $\lambda$ sufficiently large.
(If Player One keeps moving in $\tau_\lambda(g)$ until $\tau_\lambda(g)$ is finished, optimality forces Player Two to do the same.)
When $g$ and $h$ are invertible, this means that $[\tau_\lambda(g)]$ and $[h]$ are independent.
If $g$ is odd similar arguments show that $\val(\tau_\lambda(g)\oplus h)\leq \val(\tau_\lambda(g))-\val(h)$ for $\lambda$ sufficiently large.

\proclaim{Proposition 3.22}The only elements of $\overline{\sS}$ that have a non--trivial stabilizer group under the action of the $(ax+b)$--group are $\ev(x)$ for $x\geq 0$ and $[x]$ for $x\in\bR$.
\endproclaim

\demo{Proof}First suppose that an element $\alpha$ of $\overline{\sS}$ is fixed by $\varphi_{a,b}$ where $a\neq 1$.
Assume without loss of generality that $a<1$.
If $\varphi^{(n)}$ denotes the $n$--fold composition of $\varphi_{a,b}$, then $\varphi^{(n)} (x)\to b/(1-a)$, $\forall x\in\bR$.
For $\epsilon>0$ choose a graph game $g$ such that $\sI ([g],\alpha)<\epsilon$.
Then 
$$
\sI([\varphi^{(n)} (g)],\alpha)=\sI ([\varphi^{(n)} (g)],\ \varphi^{(n)} (\alpha))=a^n\sI ([g],\alpha) < a^n\epsilon.
$$
Now $[\varphi^{(n)} (g)]$ has a limit, which is either $\bold 0$ or $[b/(1-a)]$, according as $\alpha$ and $g$ are even or odd.
Therefore, either $\sI(\bold 0,\alpha)\leq 0$ or $\sI([b/(1-a)],\alpha)\leq 0$ so $\alpha=\bold 0=\ev(0)$ or $\alpha=[b/(1-a)]$.

Next suppose $\alpha$ is fixed by non--trivial $\tau_\lambda$, and without loss of generality assume $\lambda>0$.
Since $\val(\tau_\lambda(\alpha))=\val(\alpha)$, then $\alpha$ is even.
For $\epsilon>0$ choose a graph game $g$, necessarily even, such that $\sI([g],\alpha)<\epsilon$.
Then
$$
\aligned
\sI([g],[\tau_{n\lambda} (g)])&\leq \sI([g],\alpha)+\sI(\alpha,[\tau_{n\lambda} (g)])\\
& < \epsilon + \sI (\tau_{n\lambda}(\alpha), [\tau_{n\lambda} (g)])\\
& < \epsilon+\epsilon=2\epsilon.
\endaligned
$$
Now suppose $h$ is any even graph game.
Then for $n$ sufficiently large
$$
\val(g\oplus h) > -2\epsilon+\val (\tau_{n\lambda} (g)\oplus h)\geq -2\epsilon+\val(g)+\val(h).
$$
It follows that $\val(\alpha+[h])> -4\epsilon+\val(\alpha)+\val(h)$.
Since $\epsilon$ is arbitrary, we conclude that $\val(\alpha+[h])\geq\val(\alpha) +\val(h)$.
The reverse inequality follows from Proposition 3.1.
If $h$ is an odd graph game, then since $\val(g\oplus\tau_{n\lambda} (h))\leq \val (\tau_{n\lambda} (h)-\val (g))=n\lambda+\val(h)-\val(g)$ for $n$ large, it follows that $\val(\tau_{-n\lambda} (g)\oplus h)=-n\lambda+\val (g\oplus\tau_{n\lambda} (h))\leq\val(h)-\val(g)$ for $n$ large.
As above this implies that $\val(\alpha+[h])\leq \val(h)-\val(\alpha)$.
Again the reverse inequality follows from Proposition 3.1.
If $x=-\val(\alpha)$, we have shown that for all $h$, $\val ([h]+\alpha)= \val([h]+\ev(x))$; and this implies $\alpha=\ev(x)$.
\enddemo

\proclaim{Lemma 3.23}If $g$ is any graph game, then there is a number $c(g)$ such that for $y< -|g|$, $\val(g\oplus \langle y\rangle)=c(g)+y$ if $g$ is even and $\val(g\oplus \langle y\rangle)=c(g)-y$ if $g$ is odd.
If $g$ is even then $c(g)\geq |\val(g)|$, and if $g$ is odd them $c(g)\leq\val(g)$.
\endproclaim

\demo{Proof}The assumption on $y$ implies that under optimal play Player Two must never make the $y$--move if $g$ is even and Player One must never make the $y$--move if $g$ is odd.
Thus the first sentence can be proved by induction on the size of $g$, the case $g=\langle \ \rangle$ being obvious.
If $g$ is even, then $c(g)$ is the greater of $-\val(g)$ and $\max(\wt(v)-c(g\backslash v))$, where $v$ runs through all legal initial moves for $g$.
If $g$ is odd then $c(g)=\max (\wt (v)-c (g\backslash v))$.

The above already shows $c(g)\geq -\val(g)$ if $g$ is even.
To see that $c(g)\geq\val(g)$, note that Player One can save the $y$--move for last and simply play optimally in $g$.
And if $g$ is odd, Player Two can save the $y$--move for last to limit Player One to $\val(g)-y$.
\enddemo
\example{Remark}  If $f_g(x) =\val(g\oplus\langle x\rangle)$, then $f_g$ is a piecewise linear function with all slopes $\pm 1$.  It can be shown that there is (a unique) $[g_0]$ in $\sS_0$ such that $f_{g_0}=f_g$.  In the even case, if $[g_0]= \alpha + \ev(x_0)$ with $\alpha$ invertible, then $d([g],\ev(x_0))=\min\{d([g],\ev(x))\colon x\geq0\}$.  Also, $x_0 = (c(g)-\val(g))/2$.
\endexample

\proclaim{Lemma 3.24}If $g$ is an even graph game, then $c(g)=\sI (g,\langle \ \rangle)$.
\endproclaim

\demo{Proof}The main step is to show that $\val(g\oplus h)\leq c(g)+\val (h)$.
Player Two can use the following strategy:\ \ If Player One makes a move in $h$ which does not finish $h$, respond with a move optimal for $h$.
If Player One makes a $g$--move and the remnant $\tilde g$ of $g$ is odd, make a $g$--move which would be optimal in $\tilde g\oplus\langle y\rangle$ for $y<-|g|$.
In all other cases when Player Two is to move, the remnant of $g$ will be even and Player One will have previously finished $h$.
In these cases Player Two simply moves optimally in $g$.

This strategy insures that Player One's outcome from $h$ is at most $\val(h)$ and Player One's outcome from $g$ is at most $c(g)$.
One compares the $g$ portion of the play to a play of $g\oplus\langle y\rangle$ in which Player One makes the $y$--move at the moment of finishing off $h$, or at the very end of the game if Player Two finished off $h$.

Now it is easy to prove that  $ \val(g\oplus h)\geq\val(h)-c(g)$.
Player One simply makes an optimal initial move in $h$ if $h\neq\langle \ \rangle$.
This shows that $\sI(g,\langle \ \rangle)\leq c(g)$, and the reverse inequality is obvious.
\enddemo

\proclaim{Theorem 3.25}If $y< -2|g|$, then $g$ is invertible if and only if $\val(g\oplus g\oplus \langle y\rangle)=y$.
\endproclaim

\example{Remark 3.26}If $g$ is invertible, then the results of this section imply that $c(g)=\val(g)$, regardless of the parity of $g$.
Also, the main step in the proof of Lemma 3.24 is valid for general games $h$, even games without parity.
It follows that any invertible graph game is still invertible within the class of general games.
\endexample

\medskip\noindent
{\bf 3.27\ Safe play and $0-1$ games}.

In this subsection the term e--game denotes a graph game with no unbroached components and all weights equal to 0 or 1.
 The motivation is that we had some intuition  that e--games should be relatively easy to analyze, and we wanted to investigate the extent to which this is correct.
In particular Corollary 1.6 implies that if a 1--vertex is available, it is optimal to choose it.
We never had similar intuition about $0-1$ games with unbroached components.
For example, $\cyc(0,1,0,1,0,0,1,0,1,0,1,0,0,1,0,1,0,1,0,1,0)$ has obvious similarities to a game mentioned in the introduction, $\cyc(0,1,0,1,0,0,1,0,2,0,0,2,0,2,0)$.
Although these two games are not equivalent, they both have value $-1$ and both weight sequences are partitioned into what are called special slices in 4.4 below with slice weights $-2,-3,-4$.
For another example consider $\path (1,0,1,0,\ldots,1,0)\oplus \langle 1\rangle$, where the path has $n$ $1$--vertices.
There are $n+1$ vertices of weight 1, all of which are legal initial moves, and these lead to $n$ different outcomes, ranging from $n-1$ to $-(n-1)$, assuming optimal play after the initial move.
In contrast, the principle of limited badness, Lemma 1.5, implies that in an e--game  a single error can't cost a player more than two points.  Also consider the examples $g_1$, $g_2$, $h_1$, and $h_2$ from Example 1.7 with $m=1$.

If $g=(G,A)$ is an e--game, let $G_1$ be the graph obtained from $G$ by removing all the 0--vertices and the edges incident to them.
The connected components of $G_1$ will be called 1--clusters.
An equivalent game $\tilde g$ can be defined by replacing each odd 1--cluster with a single vertex of weight 1, available if and only if at least one vertex in the 1--cluster is available and adjacent to each 0--vertex adjacent to the 1--cluster, and by deleting all even 1--clusters and the edges incident to them.
New edges are added joining any two 0--vertices which are adjacent to the same even 1--cluster.
The 0--vertices in $\tilde g$ are available if and only if either they are available in $g$ or they are adjacent to an even 1--cluster which contains an available 1--vertex.
The proof that $g$ and $\tilde g$ are equivalent is by induction on the size of $g$.
Theorem 2.5 can be used to detach any available 1--vertices from $\tilde g$ and any available 1--clusters from $g$, and duplicates among the detached 1--vertices can be cancelled.
If $g$ has no available 1--vertices, then $g$ and $\tilde g$ have the same available 0--vertices, and the induction hypothesis implies that $g\backslash v\sim \tilde g\backslash v$ for such $v$.
In the game $\tilde g$ no two 1--vertices are adjacent.
We say that $\tilde g$ is reduced and is the reduction of $g$.

If $g$ is an e--game with no available $1-$vertices, we say that an available 0--vertex $v$ is a {\it safe} move if the number of 1--vertices adjacent to $v$ in the reduction of $v$ is even.
Equivalently, the number of vertices in all the 1--clusters adjacent to $v$ is even.
It is not always true that safe moves are optimal.
This probably should not be a surprise.
We proceed to identify a special class of games for which safe moves always are optimal.

We say that an e--game  is {\it simplistic} if its reduction satisfies both of the following:
\itemitem{(3.5)}All 0--vertices are available.
\itemitem{(3.6)}All non--available 1--vertices have even degree.

It is easy to check whether $g$ is simplistic without passing to the reduction.
Note that in a reduced simplistic game with no available 1--vertices the number of edges joining a 0--vertex to a 1--vertex is even.
Thus if the number of 0--vertices is odd, there is at least one safe move.
This makes it easy to prove the following by induction for simplistic games with no available 1--vertices:
\itemitem{(3.7)}If the number of 0--vertices is odd, then the value is 0 or 1, according as the number of 1--vertices is even or odd.
\itemitem{(3.8)}If the number of 0--vertices is even, then the value is 0 or $-1$, according as the number of 1--vertices is even or odd.

\noindent
The proofs of (3.7) and (3.8) show that the following strategy is optimal for both players for simplistic games:

(i)\ If there is at least one available 1--vertex, choose such a vertex.

(ii)\ If (i) does not apply and there is at least one safe move, make a safe move.

(iii)\ If neither (i) nor (ii) applies, make any legal move.

If (3.8) is applied to games of the form $g\oplus g$, $g$ simplistic, we conclude that simplistic games are invertible, with the help of Theorem 3.4.
Then if (3.7) and (3.8) are applied to $g\oplus h$, we see that any two simplistic games with no available 1--vertices and with the same parity for both 0--vertices and 1--vertices are equivalent, using Corollary 3.8.
It follows that there are exactly four equivalence classes for simplistic games with no available 1--vertices.
The simplest examples of the four classes are $\langle \ \rangle,\tes (0,1,0), \langle 0\rangle$, and $\langle 0\rangle \oplus\tes(0,1,0)$.
Of course $\tes(0,1,0)\sim\langle -1\rangle$.

We will provide four examples to show that safe play need not be optimal if only one of (3.5) and (3.6) is satisfied.
The games $g_1$ and $g_2$ satisfy (3.5) but not (3.6).
In $g_1$ there are three initial moves:\ \ One is both safe and optimal, one is safe but not optimal, and one is neither safe nor optimal.
In $g_2$ there are four initial moves, two of which are safe but not optimal and two of which are optimal but not safe.
The graph for $g_1$ is a path with weight sequence $(1,0,1,0,1,0)$ with all the 0--vertices available.
The graph for $g_2$ is obtained by adjoining two new vertices, $u$ and $w$ to a path with weight sequence $(1,0,1,0,1,0,1)$.
The vertex $w$ has weight 0 and is adjacent only to $u$.
The vertex $u$ has weight 1 and is adjacent only to $u$ and $v_4$, where $v_4$ is the middle 0--vertex.
(In $g_2$, if Player One chooses either $v_2$ or $v_6$, the outside 0--vertices, then Player Two counters with the other of $v_2,v_6$ after the exchange of 1--vertices.)
The games $g'_1$ and $g'_2$ satisfy (3.6) but not (3.5) and have corresponding initial moves, safeness status, and optimality status to $g_1$ and $g_2$.
The game $g'_1$ is obtained from $g_1$ by adjoining two non--available 0--vertices on the left, so that its graph is a path with weight sequence $(0,0,1,0,1,0,1,0)$.
The game $g'_2$ is obtained from $g_2$ by adjoining four non--available 0--vertices, two on either side.
Thus its graph is a path with weight sequence $(0,0,1,0,1,0,1,0,1,0,0)$ and with $u$ and $w$ adjoined as before.

\example{Example 3.28}We illustrate some of what is and what is not implied by the theory of metric domination and equivalence.
If the vertex $v$ is a legal initial move in $g_1$ and $w$ is a legal initial move in $g_2$, and if the $v$--move metrically dominates the $w$--move, then we have said that for any graph game $h$, the $v$--move in $g_1\oplus h$ produces at least as high an outcome for Player One as the $w$--move in $g_2\oplus h$.
If also $g_1\sim g_2$ and $w$ is an optimal move in $g_2\oplus h$, then $v$ is an optimal move in $g_1\oplus h$ and the $v$--move in $g_1\oplus h$ produces the same result as the $w$--move in $g_2\oplus h$.
But of course, if $g_1$ and $g_2$ are not equivalent, it is possible that $w$ is an optimal move in $g_2\oplus h$ and $v$ is not optimal in $g_1\oplus h$.
For a simple but unsatisfying example, take $v$ to be the 2--vertex in $g_1=\langle 2,3\rangle$, $w$ to be the 1--vertex in $g_2=\langle 1,2\rangle$, and $h=\langle 2\rangle$.
For a move interesting example, let $g_1=\cyc (0,1,0,1,1,0,1,0)$ and $g_2=\cyc (0,1,0,0,1,0)$.
Let $v$ correspond to the first 1 in $g_1$ and $w$ to the first 1 in $g_2$.
Then the $v$--move in $g_1$ and the $w$--move in $g_2$ are equivalent.
If $h=g_2$, then the $w$--move in $g_2\oplus h$ is optimal, but the $v$--move in $g_1\oplus h$ is not optimal.
And if $h=g_1$, then the $v$--move in $g_1\oplus h$ is optimal but the $w$--move in $g_2\oplus h$ is not.
For another example, let $g'_1=\cyc (0,2,0,2,2,0,2,0)$, $g'_2=\cyc(0,1,0,1,1,0,1,0)$, $v'$ correspond to a middle 2 in $g'_1$, $w'$ correspond to a middle 1 in $g'_2$, and $h=\langle 0,-1\rangle$ (cf.~Example 3.19).

Now suppose $g_1$ is a two--ended stack whose slice weight sequence, as in Lemma 2.10(i), is $(s_1,\ldots,s_m)$, and for definiteness assume $s_1\geq s_m$.
Then $g_1$ is equivalent to $g_2=\langle s_1,\ldots,s_m\rangle$.
If $v$ is the left hand move in $g_1$, then $v$ metrically dominates the $s_1$--move in $g_2$ by Corollary 2.16, and hence the $v$--move metrically dominates all moves in $g_2$.
If $h$ is any graph game, then either there is an optimal initial move from $h$ in $g_2\oplus h$ or the $s_1$--move is optimal in $g_2\oplus h$.
In the first case there is an optimal move from $h$ in $g_1\oplus h$, and in the second case the lefthand move in $g_1$ is optimal in $g_1\oplus h$.
So optimality never requires that the righthand move in $g_1$ be the first $g_1$--move in $g_1\oplus h$.
But the lefthand move in $g_1$ need not metrically dominate the righthand move, despite the fact that it metrically dominates the $s_m$--move in $g_2$, and in fact it is possible that the righthand move in $g_1$ is optimal in $g_1\oplus h$ but the lefthand move is not.
For an example take $g_1=\tes (1,2,1,2,5,2)$ and $h=\langle 5,2,0\rangle$.
Also note that the righthand move in $g_1$ is optimal in $g_1\oplus \langle 0,-1\rangle$, as are the lefthand move, the 0--move, and the $(-1)$--move.
Finally, there are two-ended stacks for which both initial moves are optimal even though $s_1 > s_m$.  For example $\tes(4,3,2,7,5)$, $\tes(4,3,0,2,7,5)$, and two examples included in 1.7.
\endexample

\vfill\eject
\noindent

{\bf \S4.\ \ Unbroached cycles and paths}

If $g$ is an unbroached cycle or path, then the best way we know to find an optimal strategy is to apply the algorithm of \S2 to $g\backslash v$ for each vertex $v$.
This gives a quadratic time algorithm which computes $\val(g)$ and determines an optimal initial move.
Once the initial move has been made, the theory of \S2 applies directly to the remaining game.
Unlike the situation for concatenations of stacks and two--ended stacks, we do not have a classification of unbroached cycles or paths up to equivalence.
However, with the help of Corollary 3.8 above and Theorem 4.3 below, we could devise a cubic time algorithm which will determine whether two given cycle games are equivalent.

\proclaim{Lemma 4.1}Let $g$ be $\cyc (a_1,\ldots,a_n)$ or $\path(a_1,\ldots,a_n)$, and let $v_1,\ldots,v_n$ be the vertices of $g$ in corresponding order.
Assume for some $i,j$ with $1\leq i< j\leq n$ that each of the following holds:

(i)\ The sequence $(a_i,a_{i+1},\ldots,a_j)$ is non--decreasing.

(ii)\ If $j<n$, then $a_j\geq a_{j+1}$; and in the cycle case, if $j=n$, then $a_j\geq a_1$.

(iii)\ If $i>1$, then $a_{i+1}\geq a_{i-1}$; and in the cycle case, if $i=1$, then $a_{i+1}\geq a_n$.

Then the $v_i$--move in $g$ is metrically dominated by the $v_j$--move.
\endproclaim

\demo{Proof}For ease of notation we divide the proof into three cases.
\enddemo

\example{Case 1:}$j=i+1$.
By Corollary 2.15, $g\backslash v_i\sim \langle a_j\rangle \oplus g\backslash \{v_i,v_j\}$.
Let $g'$ be the game obtained from $g\backslash v_j$ by changing the weight of $v_i$ to $a_j$.
Then $\sI (g', g\backslash v_j)\leq a_j-a_i$ and, again by 2.15, $g'\sim \langle a_j\rangle \oplus g\backslash \{ v_i,v_j\}$.
Therefore $\sI (g\backslash v_i,g\backslash v_j)\leq a_j-a_i$, as required.
\endexample

\example{Case 2:}$j=i+2$.
Let $g'$ be the game obtained from $g\backslash v_i$ by changing the weight of $v_{i+1}$ to $a_j$.
Then $\sI(g',g\backslash v_i)\leq a_j-a_{i+1}$ and, by Corollary 2.15, $g'\sim g\backslash \{v_i,v_{i+1}\}\oplus  a_j\rangle\sim g\backslash \{v_i,v_{i+1},v_j\}\oplus \langle a_j,a_j\rangle \sim g\backslash \{v_i,v_{i+1},v_j\}$.
Again by 2.15, $g\backslash v_j\sim \langle a_{i+1}\rangle \oplus g\backslash \{ v_{i+1},v_j\}$.
Let $g''$ be the game obtained from $g\backslash \{v_{i+1},v_j\}$ by changing the weight of $v_i$ to $a_{i+1}$.
Then by 2.15, $g''\sim g\backslash \{v_i,v_{i+1},v_j\}\oplus \langle a_{i+1}\rangle$.
Since $\sI(g'',g\backslash \{v_{i+1},v_j\}\leq a_{i+1}-a_i$, then $\sI(g\backslash v_j,\langle a_{i+1}\rangle \oplus g'')\leq a_{i+1}-a_i$.
But $\langle a_{i+1}\rangle \oplus g''\sim \langle a_{i+1}\rangle \oplus g\backslash \{v_i,v_{i+1},v_j\}\oplus\langle a_{i+1}\rangle\sim g\backslash \{v_i,v_{i+1},v_j\}$.
So the triangle inequality implies $\sI(g\backslash v_i,g\backslash v_j)\leq (a_j-a_{i+1})+(a_{i+1}-a_i)=a_j-a_i$.
\endexample

\example{Case 3:}$j> i+2$.
Let $h_i=g\backslash v_i$, $h_k=\langle a_{i+2},\ldots,a_{k+1}\rangle \oplus g\backslash \{v_i,v_{i+1},\ldots,v_k\}$ for $i< k<j$, and let $h_j=\langle a_{i+2},\ldots,a_{j-1}\rangle \oplus g\backslash \{v_i,\ldots,v_j\}$.
We claim $\sI(h_k,h_{k+1})\leq a_{k+2}-a_{k+1}$ for $i\leq k\leq j-2$.
To see this, let $h'_k$ be obtained from $h_k$ by changing the weight of $v_{k+1}$ to $a_{k+2}$ and observe that Corollary 2.15 implies $h'_k\sim \langle a_{i+2},\ldots,a_{k+1}\rangle \oplus \langle a_{k+2}\rangle \oplus g\backslash \{v_1,\ldots,v_{k+1}\}\sim h_{k+1}$.
Also note that, again by 2.15, $h_{j-1}\sim \langle a_{i+2},\ldots,a_j\rangle \oplus \langle a_j\rangle \oplus g\backslash \{v_i,\ldots,v_j\}\sim h_j$.
It follows that $\sI(g\backslash v_i, h_j)\leq\sum_i^{j-2} a_{k+2}-a_{k+1}=a_j-a_{i+1}$.
Moreover, by $j-i-1$ applications of 2.15 $g\backslash v_j\sim \langle a_{i+1},\ldots,a_{j-1}\rangle \oplus g\backslash \{v_{i+1},\ldots,v_j\}$.
Let $h''$ be obtained from $\langle a_{i+1},\ldots,a_{j-1}\rangle\oplus g\backslash \{v_{i+1},\ldots,v_j\}$ by changing the weight of $v_i$ to $a_{i+1}$.
Then $\sI(g\backslash v_j,h'')\leq a_{i+1}-a_i$ and 
$$
h''\sim \langle a_{i+1},\ldots,a_{j-1}\rangle \oplus \langle a_{i+1}\rangle\oplus g\backslash \{v_i,\ldots,v_j\}\sim h_j.
$$
It follows that $\sI(g\backslash v_j,h_j)\leq a_{i+1}-a_i$.
Thus $\sI(g\backslash v_i,g\backslash v_j)\leq (a_j-a_{i+1})+(a_{i+1}-a_i)=a_j-a_i$.
\endexample

\definition{Definition 4.2\ Plateau points}
Let $g=\path(a_1,\ldots,a_n)$ or $\cyc (a_1,\ldots,a_n)$ and let $v_1,\ldots,v_n$ be the vertices in corresponding order.
We say that $v_i$ is a {\it plateau point} of $g$ if both of the following hold:

(i)\ Either all of the vertices to the right of $v_i$ have weight $a_i$ or the first vertex to the right of $v_i$ which has a weight different from $a_i$ has weight less than $a_i$.
In the cycle case the phrase ``to the right'' is interpreted cyclically; i.e., $v_1$ is considered immediately to the right of $v_n$.

(ii)\ Same as (i) with ``right'' replaced by ``left''.

If $v_p$ is a plateau point, the largest set of consecutive vertices of equal weight which contains $v_p$ is the {\it plateau} containing $v$.
If $P_1$ and $P_2$ are two distinct plateaus such that the segment from $P_1$ to $P_2$ contains no other plateaus, we say that $P_1$ and $P_2$ are {\it neighboring}.
(In the cycle case there are two segments from $P_1$ to $P_2$, and we require only that one of these contains no other plateaus.)

If $v$ is a vertex which is not a plateau point, then there will be one or two plateaus which are considered to be {\it neighboring} to $v$.
If $P_1$ and $P_2$ are neighboring plateaus such that $v$ is in the plateau--free segment from $P_1$ to $P_2$, then each of $P_1,P_2$ is neighboring to $v$.
In the path case, if $v$ is to the left of the left--most plateau $P$, then $P$ is neighboring to $v$, and similarly if $v$ is to the right of the right--most plateau.
Finally, in the cycle case, if there is only one plateau, then it is considered neighboring to all of the non--plateau points.

There is a counterpart to the equivalence of conditions (i) and (ii) in Definition 2.9:\ \ A sequence $(a_1,\ldots,a_n)$ is weakly $U$--shaped if and only there are no interior plateaus (i.e., no plateau is entirely contained in the segment $v_2,\ldots,v_{n-1}$).
It follows in particular that the weight sequence of the plateau--free segment between two neighboring plateaus is weakly $U$--shaped.
Also, in the path case, the weight sequence to the left of the left--most plateau is non--decreasing and the weight sequence to the right of the right--most plateau is non--increasing.
And in the cycle case, if there is only one plateau $P$, then the weight sequence of the segment which starts at the right endpoint of $P$ and goes, moving to the right, around the cycle to the left endpoint of $P$ is weakly $U$--shaped.
\enddefinition

\proclaim{Theorem 4.3}Let $g=\path(a_1,\ldots,a_n)$ or $\cyc(a_1,\ldots,a_n)$, and let $v_1,\ldots,v_n$ be the vertices in corresponding order.
Then if $v$ is a vertex which is not a plateau point, there is a neighboring plateau point $w$ such that the $v$--move in $g$ is metrically dominated by the $w$--move.
Also if $w_1$ and $w_2$ are vertices which are part of the same plateau, then the $w_1$--move and the $w_2$--move are equivalent.
\endproclaim 

\demo{Proof}The first conclusion follows from Lemma 4.1 and the symmetric version of Lemma 4.1 with left and right reversed.
The only non--obvious case occurs when we have a weakly $U$--shaped sequence $(b_1,\ldots,b_m)$ with the first and last terms corresponding to plateau points and $v$ corresponding to a term $b_i$ with $1<i<m$.
(If $g$ is a cycle with only one plateau, $b_1$ and $b_m$ may correspond to the same plateau point of $g$.
This does not invalidate the reasoning.)
Let $k$ ($1\leq k < m$) be such that $(b_1,\ldots,b_k)$ is non--increasing and $(b_{k+1},\ldots,b_m)$ is non--decreasing.
Assume without loss of generality that $i>k$.
If $i\geq k+2$, Lemma 4.1 applies directly to show that the $v_i$--move is metrically dominated by the $v_m$--move (the numbering of the vertices here may have changed).
The same is true if $i=k+1$ and $b_{k+2}\geq b_k$.
If $i=k+1$ and $b_{k+2}<b_k$, then also $b_{k+1}<b_k$.
So the symmetric version of Lemma 4.1 applies to show that the $v_i$--move is metrically dominated by the $v_1$--move.

The second conclusion follows from Corollary 2.15.
\enddemo

\noindent
{\bf 4.4.\ The pizza problem}.
For brevity we will call an unbroached cycle game with non--negative weights a pizza, and we will outline our original proof that any pizza $g$ with negative value has at least fifteen vertices.
We already know that $g$ must be odd.
(See Remark 4.6 below for some additional discussion of even pizzas.)

Let $n$ be the smallest number such that there is a pizza with $n$ vertices and negative value, and let $g=\cyc(a_1,\ldots,a_n)$ be such a pizza with corresponding vertices $v_1,\ldots,v_n$.
Then: 

\noindent
(4.1) Every plateau of $g$ has length one, and for each plateau point $v_i$, we have $a_i > a_{i-1}+a_{i+1}$.
Also the segments between neighboring plateau points are $U$--shaped.

\demo{Proof}As usual, $i+1$ and $i-1$ are interpreted cyclically.
If for some $i$, $a_i\geq a_{i-1}$ and $a_i\geq a_{i+1}$, let $g'$ be the cycle of size $n-2$ obtained by condensing $v_{i-1},v_i,v_{i+1}$ to a point $w$ with weight $s=a_{i-1}-a_i+a_{i+1}$.
Then $g'$ has negative value, and we have a contradiction unless $s<0$.
To see this, note that for any move other than $w$ in $g'$ the outcome will be the same as for the corresponding move in $g$ by Corollary 2.14.
Also the outcome of the $w$--move in $g'$ will be at most the outcome of the $v_i$--move in $g$.
In fact, if Player One chooses $v_i$, Player Two must choose either $v_{i-1}$ or $v_{i+1}$, and then Player One can counter with the other of $v_{i-1}$ and $v_{i+1}$.
Since $a_i+a_{i-1}-a_{i+1}\geq s$ and $a_i+a_{i+1}-a_{i-1}\geq s$, the assertion follows.

If there were any plateaus of length at least two, then we would have the situation above with $s\geq 0$.
The fact that the segments between neighboring plateau points are $U$--shaped, rather than just weakly $U$--shaped, also follows from the previous paragraph.

Now we may also assume that $\sum_1^n a_i=1$ and that $|\val(g)|$ is maximal for pizzas of size $n$ with $\sum_1^n a_i=1$ and $\val(g)<0$.
We claim that this implies that both neighbors of any plateau point have weight 0.
To see this assume that $v_i$ is a plateau point and, say, $a_{i-1}>0$.
Then choose a sufficiently small $\epsilon>0$, and let $g'=\cyc(a_1,\ldots,a_{i-2},a_{i-1}-\epsilon,a_i-\epsilon,a_{i+1},\ldots,a_n)$.
We require that $\epsilon\leq a_{i-1}$ and $\epsilon < a_i -a_{i+1}$.  We claim that $\val(g')\leq \val(g)$.
In fact for any plateau move other than $v_i$ in $g'$, the outcome is the same as for the corresponding move in $g$, by Corollary 2.14, since $a_{i-1} -a_i +a_{i+1} = (a_{i-1}-\epsilon)-(a_i-\epsilon)+a_{i+1}$.
Also $\sI (g'\backslash v_i, g\backslash v_i)\leq\epsilon$, and this implies that the outcome of the $v_i$--move in $g'$ is at most the outcome of the $v_i$--move in $g$.
Then if $g''$ is obtained from $g'$ by multiplying all the weights by $1/(1-2\epsilon)$, we obtain a contradiction to the maximality assumption.

Now the last sentence of (4.1) implies that all the non--plateau weights are 0 and also that there are at most two 0's between any two plateaus.
Thus we can partition the weight sequence into what we will temporarily call special slices.
These are maximal length segments whose weight sequences are of the form $(0,p_1,0,\ldots,p_k,0)$ where $p_i>0$ for $i=1,\ldots,k$.
The number of special slices must be odd, since $n$ is odd.

\noindent
(4.2)\ There must be more than one special slice.

\demo{Proof}If this is false, the whole weight sequence is of the above form with $n=2k+1$.
If $w_i$ is the $i$'th plateau point, then the slice weight sequence for $g\backslash w_i$, as in Lemma 2.10(i), is just $(\ell_i,r_i)$, where $\ell_i=-\sum_1^{i-1} p_i$ and $r_i=-\sum_{i+1}^k p_i$.
Since $\ell_1=0$, $\ell_1\geq r_1$; and since $r_k=0$, $r_k\geq \ell_k$.
Therefore there is an $i$ such that $\ell_i\geq r_i$ and $\ell_{i+1}\leq r_{i+1}$.
The facts that the $w_i$---and $w_{i+1}$--moves both produce negative outcomes implies that $p_i-\ell_i+r_i<0$ and $p_{i+1}-r_{i+1}+\ell_{i+1}<0$.
Adding these two inequalities gives a contradiction.

So there are at least three special slices, and each special slice has length at least three.
The only cases with $n<15$ are $9=3+3+3$, $11=5+3+3$, $13=7+3+3$, and $13=5+5+3$.
We will leave it as an exercise for the reader to rule out these cases.
(It is a matter of checking the inconsistency of certain systems of linear inequalities.
Some cleverness can help to reduce the amount of work.)

Finally, we note that the example given in the introduction with $n=15$ has three special slices, each of length 5, and having slice--weights $-2,-3,-4$.
It is easy to check that this example has value $-1$.
\enddemo

\proclaim{Theorem 4.5}Every unbroached cycle game is invertible.
\endproclaim

\demo{Proof}This follows from Theorem 3.4$'$ (or Theorem 3.4) if we prove $\val(g\oplus g)\geq 0$.
Player One can begin by choosing a maximal weight vertex from one of the copies of $g$.
If Player Two responds with a move in the other copy of $g$, then Player One is at least even in the score and is looking at an invertible even game.
By Corollary 3.7 the outcome for Player One is non--negative.
If instead Player Two's first move is in the same copy of $g$, then Player One can choose a maximal weight vertex from the other copy of $g$.
Then after four moves Player One is at least even in the score and is looking at an invertible even game.
\enddemo

\example{Remark 4.6}This theorem combined with Corollary 3.7 gives another proof of the fact that $\val(\cyc (a_1,\ldots,a_n))\geq 0$ when $n$ is even.
Also, if $n$ is even and $\val(\cyc (a_1,\ldots,a_n))=0$, then Corollary 3.8 implies that $\cyc (a_1,\ldots,a_n)\sim\langle \ \rangle$.
When this happens, every plateau must have even length.
(If Player One chooses an initial move from a plateau of odd length, then Player One will be ahead in the score after some even sequence of moves.)
If $g=\cyc(a_1,\ldots,a_n)$ or $\path(a_1,\ldots,a_n)$ and $g'$ is the unbroached cycle or path obtained from $g$ by deleting all the vertices in a plateau $P$ of even length (this amounts to the same thing as condensing the slice consisting of $P$ and one of its neighbors),
 then $\val(g')=\val(g)$.
In fact, any move in $g'$ leads to the same outcome as the corresponding move in $g$ by Corollary 2.14.
And a move in $g$ from $P$ cannot lead to an outcome higher than $\val(g')$.
(Use Corollary 2.15 to detach the other points in $P$ and cancel any duplicates.
Then Player Two can choose the remaining vertex from $P$, and Player One is forced to choose a vertex which is an immediate neighbor $w$ of $P$.
The outcome is then the same as if Player One chose $w$ from $g'$.)
However, $g'$ need not be equivalent to $g$ (consider the $g$ of Example 3.18).
On the other hand, whenever a plateau has length greater than two, deleting two vertices from this plateau does lead to an equivalent game.

Now we see that there is an easy procedure to check whether $\val(g)=0$ when $g$ is an even unbroached cycle.
Repeatedly delete plateaus of even length from $g$ until either $g$ has become empty or there are no more even plateaus.
Some examples where $\val(g)=0$ are $\cyc(1,2,3,3,2,1)$, $\cyc(1,3,3,2,2,1)$, $\cyc(2,1,2,3,3,2,1,2)$, and $\cyc(1,3,4,4,3,2,2,1)$.
\endexample

\vfill\eject
\noindent
\example{4.7.  More on pizzas--the proportionality problem}

In this subsection we use the terminology of 4.4 and, among other things, sketch a new proof that Player  One can always obtain at least 4/9 of a pizza.  In our notation the result to be proved is that for a pizza $g$

$$
\val(g)\geq -(1/9)|g|.\tag4.3
$$

Before proceeding we make some comments.  Both of the original proofs of (4.3) in [1] and [5] are constructive and provide a good idea of what a strategy to achieve (4.3) looks like.  Our proof is totally non-constructive.  Also it would be tedious to provide all the details of our proof.  Thus, we are certainly not claiming that our proof is superior, but there are two reasons why we want to present it anyway.  One is that it shows that the methods we developed in the 1990's are powerful enough to yield a proof, and the other is that some aspects of our proof seem interesting in their own right.  In particular, we will see that our proof makes it possible to find all the extremal examples for (4.3), a problem which was discussed in [1] but only partially solved.  The paper [5] uses a concept called the heavy greens property for a segment of a pizza.  In the odd case this concept coincides with what we call a slice, and in the even case it coincides with what we call an $ev-$sequence.  However, the use made of the heavy greens property in [5] has nothing in common with the use that we make of slices and $ev-$sequences.  We believe that there is no significant overlap between our 1990's work and the work in [1] and [5].

Our proof of (4.3) has four steps.

1.  If (4.3) is false, then there is a counterexample which is partitioned into special slices.  The arguments in subsection 4.4 can be adapted with very little change to show this.

2.  If (4.3) is false, then there is a counterexample which is a $0-1$ game.  Let $g$ be a counterexample which is partitioned into special slices.  Find another counterexample $g'$ by approximating each of the non-zero weights in in $g$  by a rational number.  Then by multiplying each weight by a common denominator, we obtain a counterexample $g''$ which is partitioned into special slices and has integral weights.  Now we form a pizza $g'''$ by replacing each special slice $\underline s = (0,a_1,0,\dots,a_i,0)$ in $g''$ with a special slice $\underline s' = (0,1,\dots,1,0)$ such that $\underline s$ and $\underline s'$ have the same slice-weight (i.e., the number of 1's in $\underline s'$ is $\sum_1^i a_j$).  Then each move in $g'''$ is metrically dominated by a move in $g''$, and hence $g'''$ is also a counterexample.

3.  Let $g$ be a counterexample to (4.3) which is a $0-1$ game and which has minimal size subject to these properties.  One of the arguments in Remark 4.6 shows that it is impossible for two consecutive weights of $g$ to be 1.  (Deleting these two weights would lead to a pizza with the same value.)  And the proof of (4.1) above shows that it is impossible for three consecutive weights of $g$ to be 0.  Thus $g$ is partitioned into special slices.

Now if $g$ has $k$ $1-$vertices and $\sigma$ special slices, then it has $k+\sigma$ $0-$vertices.  Let $g'$ be obtained from $g$ by replacing each special slice $\underline s =(0,1,\dots,1,0)$ by a slice $\underline s' =(0,1,0,\dots,1,0)$ which has one fewer 1 and one fewer 0.  Note that if $\underline s = (0,1,0)$, then $\underline s' = (0)$, which is no longer a special slice.  However, it is impossible for all the special slices in $g$ to be of the form $(0,1,0)$, since then $\val(g)=1$.  (Of course $g$ is odd.)  Thus $k>\sigma$.
Finally, let $g'' =\tau_{-1}(g')$, a cycle game whose weights are all 0 and $-1$.  Each special slice $\underline s =(0,1,\dots,1,0)$ in $g$ corresponds to a slice $\underline s'' = (-1,0,-1,\dots,0,-1)$ in $g''$, and $\underline s$ and $\underline s''$ have the same slice-weight.  Clearly, every $0-$move in $g''$ is metrically dominated by a $1-$move in $g$.  Hence $\val(g'') \leq \val(g) < -k/9$.  But $\val(g'') =\val (g') - 1$, and by the minimality of $\size(g)$, $\val(g') \geq -(1/9)(k-\sigma)$.  Thus $-(1/9)(k-\sigma) -1 < -k/9$.  It follows that $\sigma < 9$.

4.  The result of step 3 allows us to reduce the proof of (4.3) to a finite problem.  For a positive integer $\sigma$ and a sequence $(x_1,\dots,x_{\sigma})$ of positive numbers we define a game $\pi(x_1,\dots,x_{\sigma})$ which is not a graph game and which has something in common with the game $\sigma (x)$ of Example 3.20.  The initial moves in $\pi(x_1,\dots,x_{\sigma})$ all have weight 0 and consist of choosing an index $i$ and a number $a$ in $[0,x_i]$.  After this initial move Player Two is faced with the graph game $\tes(-(x_i-a), -x_{i+1},-x_{i+2},\dots,-x_{i-1},-a)$, where the indices are treated cyclically (i.e., $\sigma +1$ is replaced by 1).  Then (4.3) is equivalent to the following for $\sigma =$1, 3, 5, or 7:

$$
\val(\pi(x_1,\dots,x_{\sigma})) \geq -(1/9)\sum_1^{\sigma} x_i. \tag 4.4
$$

In fact if (4.3) is false and $g$ is as in step 3, let the weights of the special slices in $g$ be $(-x_1,\dots,-x_{\sigma})$ in order.  Then each initial move in $\pi(x_1,\dots,x_{\sigma})$ is metrically dominated by a $1-$move in $g$ and $\sum_1^{\sigma} x_i =k$.  Thus $\val(\pi(x_1,\dots,x_{\sigma})) \leq \val(g)$ and (4.4) is false.  On the other hand $\pi(x_1,\dots,x_{\sigma})$ is the limit of a sequence of pizza games.  For a positive integer $n$, let $g_n = g_n(x_1,\dots,x_{\sigma})$ be the pizza obtained by replacing each $x_i$ with the weight sequence $(0,x_i/n,0,\dots,x_i/n,0)$, where there are $n$ weights equal to $x_i/n$.  Then $(g_n)$ converges to $\pi(x_1,\dots,x_{\sigma})$, and (4.4) for $\pi(x_1,\dots,x_{\sigma})$ follows from (4.3) for the $g_n$'s.
Now the verification of (4.4) for any particular value of $\sigma$ is a finite problem.  It amounts to showing that certain systems of linear inequalities are inconsistent.  Note that (4.4) is true for all values of $\sigma$.  Also the game $\pi(x_1,\dots,x_{\sigma})$ has a direct pizza-sharing interpretation.  Start with a pizza cut into pieces of sizes $(x_1,\dots,x_{\sigma})$.  Player One begins by cutting one of the pieces into two.  On each subsequent move (starting with Player Two) a player chooses one of the remaining pieces, subject to the standard selection rule, and hands it to the other player.

Now we assume that (4.3) and (4.4) are known and outline in five additional steps how to find all extremal examples for either (4.3) or (4.4).

5.  If $g=\pi(x_1,\dots,x_{\sigma})$ for $\sigma$ odd, and $m=(i,a)$ is an initial move as described above, write $g\setminus m$ for the two--ended stack which results from this initial move.  Call $m$ $i-$optimal if $\val(g\setminus m)\le \val(g\setminus m')$ for all initial moves $m'=(i,a')$.  It can be shown that for each $i$ there is an $i-$optimal move $(i,a)$ such that $(-(x_i -a),-x_{i+1},\dots,-a)$ is partitioned into slices $S_1,\dots,S_k$ with slice-weights $s_1,\dots,s_k$ such that all of the following are true:

(i) The sequence $(s_1,\dots,s_k)$ is weakly $U-$shaped.

(ii) The slice $S_1$ does not have any initial segment which is an $ev-$sequence of weight 0, and the slice $S_k$ does not have any final segment which is an $ev-$sequence of weight 0.

(iii) We have $s_1\ge s_2$ and $s_k\ge s_{k-1}$.

(iv) Either $s_1=s_k$, or $s_2=s_1>s_k$, or $s_{k-1}=s_k>s_1$.

(It can be shown that the values $s_1$ and $s_k$ are uniquely determined by conditions (i), (ii), and (iii) and that conditions (i) and (ii) together with $i-$optimality imply condition (iii).)

6.  If for some $\sigma$ there are $x_1,\dots,x_{\sigma}$ such that (4.4) is an equality, then for that $\sigma$ there are rational $x_1,\dots,x_{\sigma}$ such that (4.4) is an equality.  The game $\pi(x_1,\dots,x_{\sigma})$ could be defined also if some $x_i$'s are 0.  Then it follows from first principles that $\val(\pi(x_1,\dots,x_{\sigma}))$ is a piecewise linear function $f$ on the simplex $\{(x_1,\dots,x_{\sigma}): x_i\ge 0\medspace {\text{and}} \sum x_i =1\}$.  If $C$ is one of the cells for $f$, then $C$ is defined by linear inequalities with rational coefficients, and the minimal values of $f$ on $C$ take place on a face $C_0$.  Either one of the variables $x_i$ is identically 0 on $C_0$ or all of the variables are strictly positive on the relative interior of $C_0$.  (If $C_0$ consists of a single extreme point of $C$, then $C_0$ is equal to its relative interior.)  Since there are rational $\sigma-$tuples in the relative interior of $C_0$, there will be a rational extremal example for (4.4) in $C$ if there are any extremal examples in $C$.

7.  If $x_1,\dots,x_{\sigma}$ are positive integers, then $\val(x_1,\dots,x_{\sigma})$ is an integer.  Let $m=(i,a)$ be an optimal initial move for $g=\pi(x_1,\dots,x_{\sigma})$ for which the conditions in step 5 are satisfied.  Note that if $k>2$, then $s_2,\dots,s_{k-1}$ are integers.  If $s_1=s_k$, then $s_1$ and $s_k$ cancel out in the formula for $\val(g\setminus m)$, so $\val(g\setminus m)$ is an integer.  If, say, $s_2=s_1>s_k$, then $k>2$ and $s_1$ is an integer.  Hence $s_k$ is also an integer, and again $\val(g\setminus m)$ is an integer.  (Note that in general $\val(g)$ need not have the same parity as $\sum x_i$.  For example, $\val(\pi(1,3,5)) =0$.)

8.  If $g=\pi(x_1,\dots,x_{\sigma})$ where $x_1,\dots,x_{\sigma}$ are positive integers, and if (4.4) is an equality for $g$, let $g'$ be a $0-1$ pizza which is partitioned into special slices with slice-weights $-x_1,\dots,-x_{\sigma}$ (in order).  Then (4.3) is an equality for $g'$.

\demo{Proof}  We have $\val(g) =-r$ and $\sum x_i =9r$ for an integer $r$.  If $\val(g') > -r$, then $\val(g') \ge 2-r$ since $r$ and $9r$ have the same parity.  Consider an optimal $1-$move $v$ for $g'$, and partition the weight sequence for the two-ended stack $g'\setminus v$ into slices such that the slice-weight sequence $(s_1,\dots,s_k)$ is $U-$shaped.  Note that all of the $s_i$'s are integers.  If $s_1=s_k$, then there is a legal initial move $m$ for $g$ such that $g\setminus m$ has the slice-weight sequence $(s_1-1/2,s_2,\dots,s_{k-1},s_k-1/2)$.  If $k>2$, then $s_1\ge s_2 + 1$ and $s_k \ge s_{k-1} +1$, so in any case this new sequence is still $U-$shaped.  Thus $\val(g) \ge \val(g')-1$, a contradiction.  If, say, $s_1 > s_k$, then also $s_1 > s_2$.  Now we can choose $m$ so that $g\setminus m$ has the slice-weight sequence $(s_1 -1,s_2,\dots,s_k)$, which is weakly $U-$shaped.  Since $s_1$ is the greatest slice-weight in $g'\setminus v$ and $s_1 -1$ is the greatest slice-weight in $g\setminus m$, we see that $\val(g) \ge \val(g')$, a contradiction.
\enddemo

9.  Now we know that if there is an extremal example for (4.4) and a particular $\sigma$, then there is a $0-1$ pizza $g$ which is partitioned into $\sigma$ special slices and which is an extremal example for (4.3).  The construction in step 3 above can be applied to $g$, and it shows that $\sigma \le 9$.  Thus the determination of all extremal examples for (4.4) has been reduced to a finite problem.

We have done this finite problem without using computers.  We won't provide the details, but we will mention some tricks which shortened the work.  The main trick was to apply the methods of [1] and [5] directly to (4.4).  We actually used the methods of [5], since, even though the strategies described in [1] and [5] seem to be the same, the presentation in [5] fits better with our methods.  The difference between the $\pi-$games and pizza games mandate some changes, but there is a very close correspondence between the results of applying the methods of [5] to (4.4) and the results of [5] for (4.3).  The point is not to establish (4.4), which after all follows from (4.3), but that this method of proving (4.4) yields some equations which must be satisfied when (4.4) is an equality, and these equations are very helpful.  A second trick is to note that if the pizza $g$ in step 3 above is an extremal example and is partitioned into 9 special slices, then the pizza $g'$ of step 3 is also extremal for (4.3).  This leads to the conclusion that if $\pi(x_1,\dots,x_9)$ is an extremal example for (4.4), then so is $\pi(x_1 -\epsilon,\dots,x_9-\epsilon)$ for sufficiently small positive $\epsilon$.  This, in conjunction with the equations alluded to above, further shortens the work for $\sigma = 9$.

Our findings are that the following is the complete list of extremal examples for (4.4), up to proportionality, cyclic permutations, and reversal of direction:

(A)  $\pi(2,3,4)$

(B1)  $\pi(3,2t,3+t,9-t,12-2t)$, $0<t\le1$

(B2)  $\pi(3,2t,3+t,12-2t,9-t)$, $0<t\le1$.

Now the arguments alluded to in step 1 above show that any extremal example $g$ for (4.3) is either partitioned into special slices or is obtained from an extremal example $g'$ that is partitioned into special slices by inserting even strings of 0's into the weight sequence of $g'$.  Since these even strings of 0's will necessarily be adjacent to 0's that were already present, their insertion will not change either $\val(g')$ or $|g'|$.  If the weights of the special slices are $-x_1,\dots,-x_{\sigma}$, then one of the arguments in step 4 above shows that $\pi(x_1,\dots,x_{\sigma})$ is an extremal example for (4.4).  Now having the correct weights for the special slices is not sufficient for $g$ and $g'$ to be extremal examples for (4.3), but it is not hard to find the additional requirements by examining the outcomes of all initial plateau moves.  If $\ua =(0,a_1,\dots,a_k,0)$ is a special slice, and if $a=\sum_1^k a_i$, let $i_0$ be the smallest index such that $\sum_1^{i_0} a_i \ge a/2$.  Let $\delta(\ua) = a_{i_0} -|\sum_{i<i_0} a_i -\sum_{i>i_0} a_i|$.  Note that $\delta(\ua) \ge 0$ and $\delta(\ua)=0$ if and only if $\sum_1^{i_0} a_i = a/2$.  In case A above let $\ua$, $\ub$, and $\underline c$ be the special slices with weights $-2$, $-3$, $-4$.  Then the requirement is that $\delta(\ua)=\delta(\underline c)=0$ and $\delta(\ub) \le 1$.  In cases B1 and B2 above let $\ua$, $\ub$, $\underline c$, $\underline d$, $\underline e$ be the special slices with weights $-3$, $-2t$, $-3-t$, $-9+t$, and $-12 +2t$.  Then the requirement is that $\delta(\ub)=\delta(\underline e)=0$, $\delta(\ua)\le 3-2t$, $\delta(\underline c)\le 3-3t$, and $\delta(\underline d)\le 3-t$.

The authors of [5] asked a question near the end of the paper about the case where some weights are negative.  Under one interpretation the amount of pizza obtained by each player  is the sum of the weights of the vertices chosen by that player, and the question is whether Player One can always obtain some proportion of the total under the assumption that the sum of the weights is positive.  A counterexample is obtained with the weight sequence $(-3,2,2)$.

We prefer a different interpretation.  If $g$ is a cycle game, consider a pizza where the size of the piece corresponding to a vertex $v$ is $|\wt(v)|$.  If a player selects a vertex with negative weight, the corresponding piece of pizza is handed to the opponent.  Thus the outcome of the game, according to our usual method of scoring, is still the difference between the amounts of pizza obtained by the two players.  Under optimal play Player One gets an amount equal to $(|g| + \val(g))/2$, Player Two gets $(|g|-\val(g))/2$, and the proportion obtained by Player One is $(|g|+\val(g))/(2|g|)$.  Now if the sum of the weights is positive, Player One can obtain at least 1/4 of the pizza.  This follows from Theorem 4.9 below with $\mu = 1/2$.  And the fact that this is best possible is seen from the weight sequence $(-2,1,1)$.
\endexample

\proclaim{Lemma 4.8}Let $g$ be an odd unbroached cycle and let $m$ be the minimum weight in $g$.

(i)  If $m/|g| \leq -1/3$, then $\val(g) \geq m$.

(ii)  If $-1/3 < m/|g| < 0$, then $\val(g) \geq -(1/9)|g| + (2/3)m$.
\endproclaim
\demo{Proof}(i) Let $v_1,\dots,v_n$ be the vertices of $g$ in cyclic order, let $\wt(v_i)=a_i$, and assume $a_j = m$.  Player One uses the following strategy:  If $\sum_{j+1}^{j-1} (-1)^{i-j-1}a_i \geq 0$, initially choose $v_{j+1}$.  Otherwise take $v_{j-1}$.  (The summation is interpreted cyclically.)  On subsequent moves avoid taking $v_j$ prior to the last move of the game.  If Player Two takes $v_j$ on some move, then the outcome is at least $|m| - \sum_{i \ne j}|a_i| =2|m| - |g| \geq m$.  Otherwise, the play simply moves around the cycle ending with $v_j$ and the outcome is at least $m$.

(ii) If this is false, let $n$ be the smallest value of $\size(g)$ for which there is a counterexample.  Find a cycle $g$ of size $n$ such that $|g|=1$, $\val(g) <0$, the  minimum weight $m$ of $g$ is in $[-1/3,0]$, and $|\val(g)|/(1/9 +(2/3)|m|)$ is maximized subject to these conditions.  Furthermore assume that $g$ has the largest possible number of negative weights subject to the above.  From part (i) and the fact that (ii) is false for $g$, it follows that $m \ne -1/3$, and from (4.3), which was proved in [1] and [4], it follows that $m \ne 0$.

The minimality of $n$ implies that all plateaus in $g$ contain only one vertex.  (If $g'$ is obtained from $g$ by condensing a multiple plateau, then 4.6 implies that $\val(g')=\val(g)$.  Since $|g'| \leq |g|$, $g'$ would also be a counterexample.)  Now assume $v$ is a plateau point of $g$ with $\wt(v) >0$ and let $u$ and $w$ be its neighbors.  If, say, $\wt(u)>0$, then form $g'$ from $g$ by subtracting $\epsilon$ from $\wt(v)$ and $\wt(u)$ for sufficiently small $\epsilon > 0$.  Then $|g'| <|g|$ and $\val(g') \leq \val(g)$, a contradiction.  (See one of the arguments in 4.4.)  And if $\wt(u)=0$, the same construction also yields a contradiction, since $|g'|=|g|$ and $g'$ has more negative weights than $g$.

The conclusion from the above paragraph is that both neighbors of any plateau point of positive weight must have negative weights.  And obviously both neighbors of any other plateau point have negative weights.  Since the interval between neighboring plateau points in the weight sequence is weakly $U-$shaped, all non-plateau points have negative weight.  Let $\ell$ be the number of vertices of non-negative weight and $k$ the number of vertices of negative weight.  Since there is at least one vertex between any two neighboring plateaus, we have $k \geq \ell$, and since $g$ is odd, in fact $k> \ell$.

Now since $m$ is in the open interval $(-1/3,0)$, we are in a position to use calculus.  For sufficiently small $\lambda >0$, the game obtained by dividing all the weights of $\tau_{\lambda}(g)$ by $|\tau_{\lambda}(g)|$ is in the domain over which we took a maximum.  Let $f(\lambda)=|\val(\tau_{\lambda}(g))|/((1/9)|\tau_{\lambda}(g)| +(2/3)|m(\tau_{\lambda}(g))|)$.  Then the righthand derivative, $f'(0+)$ must be non-positive.  It is easy to see that $f(\lambda)=(|\val(g)|-\lambda)/((1/9)(1-k\lambda +\ell \lambda) +(2/3)(|m|-\lambda))$.  Thus $f'(0+)/f(0)=-1/|\val(g)| +((1/9)(k-\ell)+2/3)/(1/9 +(2/3)|m|)$.  From $f'(0+)/f(0) \leq 0$ we deduce $|\val(g)|/(1/9 +(2/3)|m|) \leq 1/((1/9)(k-\ell)+2/3)$.  Since the ratio on the left is greater than 1, we have $(1/9)(k-\ell) +2/3 < 1$, whence $k-\ell <3$.  Since $k-\ell$ is odd, we have $k-\ell =1$.

Thus there is one pair of neighboring non-negative weights which is separated by two negative weights, and all other such pairs are separated by exactly one negative weight.  So we may assume the weight sequence is $(a_1,\dots,a_n)$, where $a_i <0$ for $i$ odd and $a_i \geq 0$ for $i$ even.  In particular, for any even $j$, the sequences $(a_1,\dots,a_{j-1})$ and $(a_{j+1},\dots,a_n)$ are slices, with respective slice-weights $s_{\ell}$ and $s_r$, and $s_{\ell} = -\sum_1^{j-1}|a_i|$, and $s_r=-\sum_{j+1}^n |a_i|$.

We will now prove that $\val(g) \geq m$, a contradiction since $-1/9 +(2/3)m <m$.  Let $j$ be the largest odd number such that $\sum_1^j |a_i| \leq |m| +\sum_{j+1}^n |a_i|$.  Note that this inequality is clearly true for $j=1$ and false for $j=n$.  We claim that also $\sum_{j+2}^n |a_i| \leq |m| +\sum_1^{j+1} |a_i|$.  If this were not so, adding the two inequalities $\sum_1^{j+2} |a_i| -\sum_{j+3}^n |a_i| >|m|$ and $\sum_{j+2}^n |a_i| -\sum_1^{j+1} |a_i| >|m|$ would produce the contradiction $|a_{j+2}| >|m|$.  It follows that if Player One initially chooses $v_{j+1}$, then $|s_{\ell} - s_r| \leq |m| + \wt(v_{j+1})$.  In other words the outcome attained is at least $-|m|=m$.
\enddemo

\proclaim{Theorem 4.9}Let $g$ be an odd unbroached cycle game, and let $\mu$ be a number in the interval $[0,1]$.

(i) If $\mu \geq 1/3$ and $\wt(v) \geq -\mu |g|$ for all vertices $v$ of $g$, then $\val(g) \geq -\mu |g|$.

(ii) If $0 \leq \mu \leq 1/3$ and $\wt(v) \geq -\mu |g|$ for all vertices $v$ of $g$, then $\val(g) \geq -(1/9)|g| -(2/3)\mu |g|$.
\endproclaim

\demo{Proof}Let $f(\mu)=1/9 + (2/3)\mu$ for $\mu \leq 1/3$ and $f(\mu) = \mu$ for $\mu \geq 1/3$.  Note that $f$ is continuous and monotone increasing.  Then part (i) follows from both parts of the lemma combined with the monotonicity of $f$, and part (ii) for $\mu >0$ follows from part(ii) of the lemma and the monotonicity.  Of course, part (ii) for $\mu =0$ is (4.3), which was proved in [1] and [5], and which was used in proving the lemma.
\enddemo

\example{4.10.  Examples and remarks}The estimates in the theorem are sharp.  For $\mu \geq 1/3$ a simple example is $\cyc(-2\mu,1-\mu,1-\mu)$.  Also for $\mu \geq 1/3$ it can be shown that the most general extremal example is $\cyc(a_0,a_1,\dots,a_{n-1})$ where $a_0 =-\mu|g|$, $\val(\path(a_1,\dots,a_{n-1}))=0$, and $a_i\geq -\mu|g|,\forall i$.  For $\mu = 1/3$ the proof is complicated by a somewhat tricky phenomenon:  For example $\cyc(-2,-1,-2,1,0)$ is an extremal example despite the fact that $\path(-1,-2,1,0)$ has positive value.  However, $\path(1,0,-2,-1)$ does have value 0.

For $\mu<1/3$ any extremal example must have size at least 15, since otherwise $\val(g) \geq m$.  Let $g_k = \cyc(0,1,-(k-2),1,0,0,k-1,0,2,0,0,2,-(k-2),2,0)$ for $2\leq k \leq 3$ and $g_k = \cyc(0,1,-(k-2),1,0,0,2,-(k-3),2,0,0,2,-(k-2),2,0)$ for $k>3$.  Thus $g_k$ is partitioned into slices of weights $-k$, $-(k+1)$, and $-(k+2)$.  Then $\val(g_k) = -(k-1)$.  This is easily proved by applying the theory in $\S 2$ to all possible initial plateau moves.  Thus $g_k$ is an extremal example for $\mu=(k-2)/(3k+3)$, and the whole interval $[0,1/3)$ is covered.  Of course $g_2$ is the example mentioned in the introduction.  Also note that if we divide the weights of $g_k$ by $|g_k|$ and let $k \to \infty$, we get a game with three weights equal to -(1/3), a different extremal example for $\mu = 1/3$ from  either of the ones given above.

It is amusing to consider odd cycle games $g$ whose weights are all $\pm 1$ and to ask what is the minimum of $\val(g)$ for a given value of $n=\size(g)=|g|$.  Of course, this is equivalent to asking what is the minimum of $\val(g')$ for $0-1$ pizzas $g'$ with $n$ vertices.  The answer to the first question, for $n\ge 3$, is $-2\lfloor (n-3)/18\rfloor -1$.  One way to see this is to apply Theorem 4.9(ii) for $\mu = 1/n$ and note that $-2\lfloor (n-3)/18\rfloor -1$ is the smallest odd integer which is at least $-(1/9)n -(2/3)(1/n)n$.  This shows that the minimum is at least $-2\lfloor (n-3)/18\rfloor -1$.  It is easy to get extremal examples if $n<21$.  If $\lfloor (n-3)/18\rfloor =k>0$, start with a $0-1$ pizza $g'$ which is partitioned into 3 special slices with weights $-2k$, $-3k$, and $-4k$.  Then $\size(g') =18k+3$ and $\val(g')=-k$.  Obtain $g''$ with the same value by inserting a string of $n-18k-3$ consecutive 0's, and obtain $g$ by changing each weight $w$ in $g''$ to $2w-1$.  If $n=18k+3$, this is the only way to get extremal examples, but otherwise there are additional extremal examples.

It might be more in the spirit of the question asked in [5] to consider a different question from the one answered in Theorem 4.9.  Let $N=\sum_{\wt(v)<0}|\wt(v)|$.  For $\nu$ in $[0,1]$, what is the infimum of $\val(g)/|g|$ for odd cycles with $N\leq \nu |g|$?  Our interpretation of the actual question in [5] corresponds to $\nu = 1/2$.  Since $m \geq -N$, it is obvious that the infimum must be at least $-f(\nu)$ for $f$ as in the proof of 4.9.  And if $\nu \geq 1/3$, this is the correct answer, since there are extremal examples for Theorem 4.9(i) with only one negative weight.  It can be shown that extremal examples for Theorem 4.9 with $0<\mu <1/3$ must have at least two negative weights.  What we can say for the other values of $\nu$ is that the infimum is between $-(1/9) -(2/3)\nu$ and $-(1/9) -(1/3)\nu$ for $0<\nu < 1/6$ and between $-(1/9) -(2/3)\nu$ and $- \nu$ for $1/6 \leq \nu <1/3$,

There is little to be said about even cycles in this sort of context.  We know that all even cycle games have non-negative value, and Remark 4.6 above gives a reasonable idea what the value 0 cases look like.  It may be somewhat interesting to note that if $g$ is an even $0-1$ cycle, then $\val(g) = 0$ if and only if all the $1-$clusters in the weight sequence are even, and $\val(g) =1$ if and only if there is exactly one odd $1-$cluster. Also, if one of the weights has magnitude greater than $(1/2)|g|$, then the value must be positive.  The related inequalities are too obvious to be worth stating formally.

On the other hand, $\val(\pi(x_1,\dots,x_{\sigma})) \ge (1/\sigma)\sum x_i >0$ for even $\sigma$.  Equality occurs if and only if all the $x_i$'s are equal.

\demo{Proof}  Perform a cyclic permutation to achieve $\sum(-1)^{i-1}x_i =\delta \ge 0$ and $x_2\ge x_j$ for all even values of $j$.  Player One begins with the move $m=(2,x_2/2)$.  Player Two must choose one of the vertices of weight $-x_2/2$ and Player One follows with the other.  Thereafter, Player One follows a color strategy and chooses only vertices with even indices.  Thus the outcome for Player One is $\sum_1^k x_{2i-1}-\sum_2^k x_{2i} = x_2 + \delta$ if $\sigma =2k$. Since $\sum_1^k x_{2i} = ({\sum_1^{\sigma} x_i}-\delta)/2$, $x_2 \ge (({\sum_1^{\sigma} x_i}-\delta)/2)/(\sigma/2) =(1/\sigma)\sum_1^{\sigma} x_i -(\delta/\sigma)$.  Thus $\val(\pi(x_1,\dots,x_{\sigma})) \ge (1/\sigma)\sum_1^{\sigma} x_i + (1-(1/\sigma))\delta \ge (1/\sigma)\sum_1^{\sigma} x_i$.  If equality occurs then $\delta = 0$ and all even $x_i$'s are equal.  But if $\delta = 0$, we can apply the same strategy with even and odd reversed to conclude that also all odd $x_i$'s are equal.
\enddemo
\endexample

\proclaim{Theorem 4.11}If $g$ is an unbroached path game, then the following are equivalent:

(i) For each vertex $v$ of $g$ there is a vertex $w$ such that the $v$--move is metrically dominated by the $w$--move and $g\backslash w$ is invertible.

(ii) For each plateau point $v$ there is a plateau point $w$ such that the $v$--move is metrically dominated by the $w$--move and $g\backslash w$ is invertible.

(iii) The game $g$ is invertible.
\endproclaim

\demo{Proof}It follows from Theorem 4.3 that (ii) implies (i).
To see that (i) implies (iii), we use Theorem 3.4$'$.
It is enough to prove $\val(g\oplus g)\geq 0$.
Player One can choose a vertex $v$ from one of the copies of $g$ such that $\wt(v)$ is maximal and $g\backslash v$ is invertible.
If Player Two responds in the other copy of $g$, then (i) implies that Two can do no better than to choose a vertex $w$ such that $g\backslash w$ is invertible.
Then, as in the proof of Theorem 4.5, the outcome for One will be non--negative.
If instead Player Two responds in the same copy of $g$, then Player One can choose the $v$ in the other copy.
Again the arguments from 4.5 apply.

To prove that (iii) implies (ii) we consider $h=g\oplus g\oplus\langle y\rangle$ where $y<<0$.
Then optimality requires that Player Two never makes the $y$--move, and invertibility of $g$ implies that $\val(h)= y$.
Suppose Player One chooses a plateau point $v$ of $g$ such that $v$ is not metrically dominated by a plateau point $w$ with $g\backslash w$ invertible and such that $\wt(v)$ is maximal for such choices of $v$.
If Player Two responds in the same copy of $g$, then One can choose a maximal weight vertex from the other copy of $g$.
Regardless of Two's reply, One will be at least even in the score after four moves and will be looking at a game of the form $k\oplus \ev(x)\oplus\langle y\rangle$, where $k$ is an invertible even game and $x>0$.
The $\ev(x)$ can be ignored for the further strategy, and One just continues with optimal play in $k$.
So the outcome for One is at least $x+y$, a contradiction.

So we may assume Player Two's first move is a vertex $w$ from the other copy of $g$, and by Theorem 4.3 we may assume $w$ is a plateau point.
Let $g\backslash v\sim g'\oplus \ev(x')$ where $g'$ is invertible and $x'>0$ and $g\backslash w\sim g''\oplus \ev(x'')$ where $g''$ is invertible and $x''\geq 0$.
Then Player One's outcome is $\sI(g',g'')+x'+x''+y+\wt(v)-\wt(w)$, by Corollary 3.8.
It follows that $x'+x''+\sI(g',g'')\leq \wt(w)-\wt(v)$.
Since $\sI(g\backslash v, g\backslash w)=\sI (g',g'')+ |x'-x''|$ by Proposition 3.12, we conclude that the $w$--move metrically dominates the $v$--move.
Since also $\wt(w) > \wt(v)$, the maximality hypothesis for $\wt(v)$ implies that the $w$--move is metrically dominated by a plateau move $u$ with $g\backslash u$ invertible.  Then the $u$--move metrically dominates the $v$--move, a contradiction.  
\enddemo

\proclaim{Proposition 4.12}Let $g=\path (a_1,\ldots,a_n)$ and $h=\cyc (a_1,\ldots,a_n)$.

(i)\ If $n$ is odd, then $\val(g)\geq\val(h)$; and if $n$ is even, then $\val(g)\leq \val(h)$.

(ii)\ The game $g$ is equivalent to the game $h$ if and only if $g$ is invertible.

(iii)\ If $(a_1,\ldots,a_n)$ is a slice, then $g\sim h$.
Moreover, if $n$ is odd some cyclic permutation of $(a_1,\ldots,a_n)$ is a slice.
Thus for odd $n$, $\val(h)=\min (\val (\path (\ub))$, where $\ub$ runs through all cyclic permutations of $\ua$.

(iii$'$)\ If $n$ is even, then $\val(h)=\max (\val (\path (\ub))$, where $\ub$ runs through all cyclic permutations of $\ua$.
\endproclaim

\demo{Proof}(i)\ We compare $g\backslash v$ and $h\backslash v$.
By Lemma 2.10 and Theorem 2.11 $g\backslash v\sim \st(s'_1,\ldots,s'_m)\oplus \st(\ux')\oplus \st (s''_1,\ldots,s''_p)\oplus\st(\ux'')$, and $h\backslash v\sim\tes (s'_1,\ldots,s'_m,\ux',(\ux'')^{-1}$, $s''_p,\ldots,s''_1)$, where $(s'_1,\ldots,s'_m)$ and $(s''_1,\ldots,s''_p)$ are decreasing, $\ux'$ and $\ux''$ are $\ev$--sequences, and $(\ux'')^{-1}$ denotes the reverse of $\ux''$.
(Possibly $m$ or $p$ is 0, and possibly $\ux'$ or $\ux''$ is missing.)
Let $x'$ be the weight of $\ux'$ and $x''$ be the weight of $\ux''$.
Let $\uy'$ be obtained from $\ux'$ by subtracting $x'$ from the last term of $\ux'$, and let $\uy''$ be obtained from $\ux''$ in a similar way.
Let $k=\tes (s'_1,\ldots,s'_m,\uy',(y'')^{-1},s''_p,\ldots,s''_1)$.
Then $\sI (h\backslash v,k)\leq x'+x''$ and $(\uy',(\uy'')^{-1})$ is an $\ev$--sequence of weight 0.
It follows that $(\uy',(\uy'')^{-1})=(S,S')$ where $S$ and $S'$ are slices of the same weight $s$.
Regardless of the value of $s$, $k\sim \langle s'_1,\ldots,s'_m,s''_p,\ldots,s_1^k\rangle$.
(If, say, $s\geq s'_m$, then $(s'_m,s,s)$ is a slice which can be condensed to $s'_m$.
And if $s< s'_m$, then the sequence $(s'_1,\ldots,s'_m,s,s,s''_p,\ldots,s''_1)$ is weakly $\cU$--shaped.)

So if $n$ is odd, then $\val(g\backslash v)=\val (\langle s'_1,\ldots,s'_m,s''_1,\ldots,s''_p\rangle)-x'-x''$ and $\val(h\backslash v)\geq \val(\langle s'_1,\ldots,s'_m,s''_1,\ldots,s''_0\rangle)-x'-x''$.
The fact that $\val(h\backslash v)\geq \val(g\backslash v)$, for all $v$ implies $\val(h)\leq \val(g)$.
If $n$ is even, a similar argument applies.

(ii)\ Suppose $v$ and $w$ are vertices such that $g\backslash w$ is invertible and the $w$--move in $g$ metrically dominates the $v$--move in $g$.
We claim that also the $w$--move in $h$ metrically dominates the $v$ move in $h$.
In the proof of part (i) we saw that $g\backslash w\sim h\backslash w$.
We also saw that $g\backslash v\sim g'\oplus\ev (x)$ where $g'$ is invertible and $\sI (g', h\backslash v)\leq x$.
By 3.12 and the given metric domination, $\sI(g',g\backslash w)+x\leq \wt (w)-\wt(v)$, and this implies the claim.

Now it is obvious that $g\sim h$ implies $g$ invertible.
If $g$ is invertible, by Corollary 3.8 it is sufficient to prove $\val(g\oplus h)\leq 0$.
But by the above paragraph and Theorem 4.11, we may assume that Player One's initial move is a vertex $v$, taken from either $g$ or $h$, such that $g\backslash v$ is invertible.
Player Two can then choose the ``same'' vertex $v$ from the other game.
Since $g\backslash v\sim h\backslash v$ for such a vertex, the result follows.

(iii)\ If $(a_1,\ldots,a_n)$ is a slice, then the weights called $x$ and $x'$ in the proof of (i) must both be 0.
Thus $g\backslash v$ is invertible for all $v$, $g$ is invertible, and $g\sim h$.
To see that some cyclic permutation is a slice, we use induction on $n$, the case $n=1$ being obvious.
If $n\geq 3$, perform a cyclic permutation so that $a_2$ is the maximal weight.
Then apply the induction hypothesis to find a cyclic permutation of $(s,a_4,\ldots,a_n)$ which is a slice, where $s=a_1-a_2+a_3$.
Then it follows from Lemma 2.8 that the corresponding cyclic permutation of $(a_1,\ldots,a_n)$ is a slice.

(iii$'$)\ Let the vertex $v$ be an optimal initial move for $h$.
Then $h\backslash v=\tes (S_1,\ldots,S_m)$ where $S_1,\ldots,S_m$ are slices whose weight sequence $(s_1,\ldots,s_m)$ is $U$--shaped.
Choose $j$ such that $s_j$ is minimal, and permute $(a_1,\ldots,a_n)$ cyclicly in such a way as to obtain $(S_{j+1},\ldots,S_m,\wt(v),S_1,\ldots,S_j)$.
Then the $v$--move in the permuted path game produces an outcome equal to $\val(h)$.
\enddemo
\example{Remark 4.13}  With the notation of the proposition, it is an easy exercise to show that $\path(\ub)$ is invertible for all cyclic permutations of $\ua$ if and only if either all the $a_i$'s are the same or all but one are the same with the odd one being higher.  It is also possible, though harder, to find necessary and sufficient conditions for $\val(\path(\ub))$ to be the same for all $\ub$.  In the even case the condition is that the color strategy be optimal for the cycle game $h$.  In the odd case it is that that $h$ have only one plateau, which is odd, and that $h\backslash v\sim \langle \ \ \rangle$ for a plateau point $v$.  When this happens, the greedy strategy is optimal for both players throughout the play of $g$ or $h$.  A key fact used in the proofs for both parities is that if, in the proof of 4.12(i), both $x'>0$ and $x''>0$, then the inequality proved there is strict.
\endexample

\example{Examples}If $g=\path(1,2,3)$, then $g$ is invertible despite the fact that $(1,2,3)$ is not a slice.
(Of course this implies $g\sim\path(2,3,1)$.)
But $\path(3,1,2)$ is not invertible and $\val(\path(3,1,2))=4>2=\val(\path(1,2,3))$.
If $g=\path(1,0,1,1,0)$, then $g$ is not invertible, but $\val(g)=1=\val(\cyc (1,0,1,1,0))$.
(Of the five cyclic permutations of this sequence, two are slices, two give non--invertible path games of value 1, and one gives a path game of value 3.)
If $g=\path(0,1,0,2)$, then $g$ is non--invertible, but $\val(g)=3=\val (\cyc (0,1,0,2))$.
In this case, all of the cyclic permutations give non--invertible path games of value 3.
The only truly different cyclic permutation is $(0,2,0,1)$ and $g$ is not equivalent to $\path(0,2,0,1)$.
(Note that $\val(g\oplus \langle -1,-2\rangle)=4$ and $\val(\path(0,2,0,1)\oplus \langle -1,-2\rangle)=2$.)
If $g=\path (1,2,4,3)$, then $g$ is invertible and has value 2.
One of the non--trivial cyclic permutations is also invertible, and the other two give path games of value 0.
If $g=\path(1,2,3,4,5)$, then $g$ is invertible and has value 3.
One of the cyclic permutations, namely $(2,3,4,5,1)$ is a slice.
The other three cyclic permutations give non--invertible path games of values 5, 3, and 7.
If $g_1=\path(0,1,0,2,3,2)$ and $g_2=\path (0,1,0,0,2,3,2)$, then each $g_i$ is invertible despite the face that $g_i$ has a plateau point $v$ with $g_i\backslash v$ not invertible.
The game $g_1$ has value 4, one of its cyclic permutations is also invertible, and the other four cyclic permutations have value only 2.
The game $g_2$ has value 2, one of the cyclic permutations has a weight sequence which is a slice, and the other five cyclic permutations have value 4.
\endexample

\bigskip\noindent
{\bf \S5.\ Some variants of graph games}

\noindent
{\bf 5.1\ Go--style passing}.

Under these rules, which will be called $p$--rules, a player may pass at any turn instead of removing a vertex.
The move then passes to the opponent, and the game ends if there are two consecutive passes or when all vertices have been removed.
For a graph game board $g$, the $p$--rules game played on $g$ will be denoted by $\Gamma_p(g)$, and $\val (\Gamma_p(g))$ will also be denoted by $\val_p(g)$.
Then $\val_p(g)=\max (0,\max (\wt (v)-\val_p (g\backslash v)))$, where in the second max $v$ ranges through all legal initial moves for $g$.
Whenever it is optimal for a player to pass, it is optimal for the opponent to respond with a pass.
If we specified that $n$ consecutive passes end the game, $\val_p (g)$ would be the same so long as $n>0$.

Note that $\Gamma_p (g_1\oplus g_2)\neq \Gamma_p (g_1)\vee \Gamma_p (g_2)$.
In the latter game, the concatenation of two $p$--rules graph games, it is possible, for example, that one of $g_1,g_2$ can be passed out while the other is still active.
However, it can be shown that $\val_p (g_1\oplus\ldots\oplus g_n)=\val (\Gamma_p (g_1)\vee \ldots\vee\Gamma_p (g_n))$.

We say that $g_1$ is $p$--equivalent to $g_2$, $g_1\underset p\to\sim g_2$, if $\val_p (g_1\oplus h)=\val_p (g_2\oplus h)$ for all graph game boards $h$, and define $\sI_p (g_1,g_2)=\sup|\val_p (g_1\oplus h)-\val_p (g_2\oplus h)|$.
Under $p$--rules, it is never optimal to select a vertex of negative weight, and it can't be optimal to select a vertex of weight 0 unless it is also optimal to pass.
If all the vertices of $g$ have non--positive weight, then $g\underset p\to\sim \langle \ \rangle$, regardless of the parity of $g$.
The parity of  $\Gamma_p(g)$ is not predetermined.

Despite the lack of parity, a number of positive results hold.
For all $g$, $g\oplus g\underset p\to\sim \langle\ \rangle$.
And the inequalities (3.2) and (3.3), with $\val$ replaced by $\val_p$, hold for all $g$ and $h$.
Corollary 2.4, Theorem 2.5, and Theorem 2.6 all hold with $\sim$ replaced by $\underset p\to\sim$, but now $\st(a_1,\ldots,a_n)\underset p\to\sim \langle\ \rangle$ if $(a_1,\ldots,a_n)$ is an $\ev$--sequence.
Thus we still have a complete analysis of concatenations of stacks and two--ended stacks under $p$--rules.
However, $\val_p (\langle s_1,\ldots,s_m\rangle)$ may be different from $\val(\langle s_1,\ldots,s_m\rangle)$ if some $s_i$'s are negative.
On the other hand, $g_1\underset p\to\sim g_2$ does not imply $\tau_\lambda (g_1)\underset p\to\sim\ \tau_\lambda(g_2)$, and therefore we have no translation operation on the group of $p$--equivalence classes.
(See the next subsection for more on this.)

\medskip\noindent
{\bf 5.2.\ Asymptotic rules}.

Under these rules, which will be called $a$--rules, a player may pass only if the number of remaining vertices is even, and two consecutive passes still end the game.
So $\val_a(g)=\max (0,\max (\wt (v)-\val_a (g\backslash v)))$ if $g$ is even and $\val_a (g)=\max (\wt (v)-\val_a (g\backslash v))$ if $g$ is odd.
It is not literally true that the parity of $\Gamma_a(g)$, the $a$--rules version of $g$, is pre--determined, since it is not legally required that the number of passes be even.
However, for determining $\val_a(g)$, we can assume that any pass is immediately followed by another pass; and hence the parity of $g$ is relevant.

The value of $\Gamma_a (g_1)\vee \Gamma_a (g_2)$ need not be the same as $\val_a(g_1\oplus g_2)$.
(For example take $g_1=\langle 0,-2\rangle$ and $g_2=\langle 1\rangle$.)
This is not a problem.
The equivalence relation, distance function, and addition operation on equivalence classes are all defined via $\oplus$, not via $\vee$.

It is not hard to see that for $\lambda$ sufficiently large, $\val_p(\tau_\lambda(g))=\val_a (g)$ if $g$ is even and $\val_p (\tau_\lambda (g))=\lambda+\val_a (g)$ if $g$ is odd.
This explains the word ``asymptotic'' and enables us to deduce results for $a$--rules from results for $p$--rules.
In particular $g\oplus g\ \underset a\to\sim\ \langle \ \rangle$ for all $g$, and the inequalities $(3.2)$ and $(3.3)$, with $\val$ replaced by $\val_a$, hold under the same parity assumptions as for ordinary invertible graph games.
Also if all the vertices have weight 0, then $g\ \underset a\to\sim\ \langle\ \rangle$ if even and $g\ \underset a\to\sim\ \langle 0\rangle$ if odd.

It is easy to prove by induction that $\val_a(g)\geq\val(g)$ if $g$ is even, and $\val_a(g)\leq\val(g)$ if $g$ is odd.
Also the analogues of Lemmas 3.23 and 3.24 hold for $a$--rules, with $\val_a(g)$ playing the role of $c(g)$.
Therefore if $g$ is an even graph game board, then
$$
\val(g)\leq\val_a(g)=\sI_a(g,\langle\ \rangle)\leq \sI(g,\langle \ \rangle);\tag5.1
$$
and if $g$ is an odd graph game board, then
$$
c(g)\leq \val_a (g)\leq \val(g).\tag5.2
$$
In particular $\val_a(g)=\val(g)$ if $\Gamma(g)$ is invertible (using ordinary rules).
Applying (5.1) to $g\oplus h$, where $g$ and $h$ have the same parity, we obtain $\sI_a(g,h)\leq\sI (g,h)$ if at least one of $\Gamma(g)$, $\Gamma(h)$ is invertible.
It follows that there is an isometric homomorphism $\Phi\colon \sG\to \sG_a$ where $\sG_a$ is the group of $a$--equivalence classes.

Also, $\val_a(g) \geq \val_p(g)$ if $g$ is even, and $\val_a(g) \leq \val_p(g)$ if $g$ is odd.
Thus $\sI_p(g,h)=\val_p (g\oplus h)\leq\val_a (g\oplus k)=\sI_a (g,h)$ when $g$ and $h$ have the same parity, whence there is a contractive surjective homomorphism $\Psi\colon \sG_a\to \sG_p$.

It would be nice if the results of the above paragraph could be improved.
Here are some questions motivated only by hope without any reason for believing the answers are affirmative:
$$
\text{Is }\sI_a(g,h)\leq \sI(g,h)\text{ even when }\Gamma(g)\text{ and } \Gamma(h)\text{ are non-invertible?}\tag5.3
$$
$$
\text{Does }g\sim h\text{ imply }g\ \underset a\to\sim\ h\text{?}\tag5.3$'$
$$
$$
\text{Is }\Phi\text{ surjective?}\tag5.4
$$
$$\aligned
\text{For every }(a_1,\ldots,a_m),\text{ does there exist }(b_1,\ldots,b_n)\\
\text{ such that } \path (a_1,\ldots,a_m)\ \underset a\to\sim\ \cyc (b_1,\ldots,b_n)\text{?}\endaligned\tag5.4$'$
$$
$$
\text{Is }\val_p(g)\text{ equal to} \inf\{\val_a(h)\colon h\text{ is even and }\Psi([h]_a)=[g]_p\}\text{?}\tag5.5
$$
If (5.3$'$) and (5.4) are affirmative, then the map $\Phi$ would extend to $\sS$ and would amount to a retraction from $\sS$ onto $\sG$.
If (5.5) is affirmative, then not only is $\sG_p$ a quotient group of $\sG_a$, but $\sI_p$ would be the quotient metric induced from $\sI_a$.
It can be shown that (5.4$'$) is affirmative if $m\leq 4$.
For example, $\path(1,0,1,0)\ \underset a\to\sim\ \cyc(1/2,1,1/2,0,1,0)$.

Our hope with (5.3) and (5.4) is to learn more about graph games under ordinary rules, and partial results (such as (5.4$'$)) would be welcome.
However, if one's goal were just to construct a large well--behaved group of equivalence classes, one might just make $\sG_a$ the primary object of study.

\noindent
{\bf 5.3.\ $s$--rules}.
Under $s$--rules passing is allowed but two consecutive passes are not allowed.
If one player passes, then the other is required to choose a vertex on the next move.
In particular Player One can force Player Two to take all the vertices, and this strategy is sometimes optimal.
If $\valsm(g)=\max (\wt (v)-\val_s (g\backslash v))$, where $v$ runs through all legal initial vertex moves in $g$, then $\val_s(g)=|\valsm(g)|$.
There is no pre--determined parity under $s$--rules.
As with $a$--rules, equivalence, distance, and addition are defined in terms of $\oplus$, not $\vee$.

It seems difficult to prove strong positive results about $s$--rules.
The game $\Gamma_s(g)$ is $s$--equivalent to $\ \langle\ \rangle$ if and only if all weights in $g$ are 0.  In particular there are no non--trivial invertible equivalence classes under $s$--rules.  
The value of $g=\langle a_1,\ldots,a_k\rangle$ under $s$--rules can be calculated as follows:\ \ Let $N=\sum_{a_i<0} |a_i|$ and $P=\sum_{a_i>0} a_i$.
If $N>P$, then $\val_s(g)=N-P$.
If $P\geq N$, let $(p_1,p_2,\ldots,p_k)$ be a listing of the positive weights in increasing order, and let $j$ be the smallest index such that $P'=\sum_1^j p_i\geq N$.
Then $\val_s(g)=p_k-p_{k-1}+\ldots+ (-1)^{k-j+1}p_{j+1}+(-1)^{k-j} (P'-N)$.
This formula can be used to construct some strange--looking examples.

\noindent
{\bf 5.4\ Move menus}.
Define an abstract move as an ordered pair $(a,g)$ where $a\in\bR$ and $g$ is a graph game board.
Let $M$ be a finite set of abstract moves, all of the same parity.
In the game $\Gamma(M)$ based on $M$, the initial move is the choice of an element $(a,g)$ of $M$.
Then $a$ becomes Player One's first score, $-a$ is Player Two's first score, and Player Two is then faced with the graph game $g$.
Thus $\val(\Gamma(M))=\max (a-\val(g))$, where $(a,g)$ runs through $M$.

\noindent
{\bf 5.5\ Even vertices}.
We can generalize graph games by regarding some of the vertices as even.
If a player chooses an even vertex in the course of play, then that player takes another turn.
If one wanted to make a serious study of graph games with even vertices, one should modify our previous notation $\ev(x)$.

\bye